\documentclass[11pt]{amsart}
\usepackage{amsmath,amssymb,amsthm,amsxtra,graphics,enumerate} 
\usepackage{stmaryrd}
\usepackage{times}
\usepackage[margin=1.3in]{geometry}
\usepackage[linktocpage=true, hidelinks, colorlinks=false, allcolors=blue]{hyperref}
\usepackage{comment}
\usepackage{scalerel}
\usepackage{mathrsfs}

\newtheorem{thm}{Theorem}[section]

\newtheorem{cor}[thm]{Corollary}

\newtheorem{lem}[thm]{Lemma}
\newtheorem{prop}[thm]{Proposition}

\newtheorem{question}[thm]{Question}
\theoremstyle{definition}
\newtheorem{defn}[thm]{Definition}
\newtheorem{rem}[thm]{Remark}

\newcommand\meet{\wedge}

\newcommand\Seq[1]{( #1 )}
\newcommand\one{\mathbf{1}}
\newcommand\St{{}^{\mathrm{st}}}

\newcommand\Ecal{\mathcal{E}}
\newcommand\Fcal{\mathcal{F}}
\newcommand\Gcal{\mathcal{G}}
\newcommand\Hcal{\mathcal{H}}

\newcommand\Ncal{\mathcal{N}}

\newcommand\p{\mathbb{P}}
\newcommand\q{\mathbb{Q}}

\newcommand\forces{\Vdash}

\newcommand\Rscr{\mathscr{R}}

\newcommand\Rbb{\mathbb{R}}
\newcommand\Sbb{\mathbb{S}}

\renewcommand\emptyset{\varnothing}

\newcommand\dom{\operatorname{dom}}

\begin{document}

\title[Five element basis]{A virtual five element basis for the \\ uncountable linear orders}

\author{John Krueger and Justin Tatch Moore}

\address{John Krueger, Department of Mathematics, 
	University of North Texas, Denton, TX, USA}
\email{john.krueger@unt.edu}

\address{Justin Tatch Moore, Department of Mathematics, Cornell University, 
Ithaca, NY, USA}
\email{justin@math.cornell.edu}

%\subjclass{}

%\keywords{}

%\date{September 29, 2025}

\begin{abstract}
    We prove that for every Aronzsajn line $A$ and every Countryman line $C$, there is a proper forcing extension in which $A$ contains an
    isomorphic copy of either $C$ or its converse $C^*$.
    As a corollary, we obtain answers
    to several related questions asked by the second author in the literature:
    if there is an inaccessible cardinal,
    then there is a proper forcing extension in which the uncountable linear orders have a five element basis;
    $\textsf{BPFA}$ implies
    the existence of a five element basis
    for the uncountable linear orders;
    $\textsf{BPFA}$ is equiconsistent with
    the conjunction of $\textsf{BPFA}$ and
    Aronszajn tree saturation.
    These results are derived from
    new preservation results
    concerning subtrees of Aronszajn
    trees, proper forcings, and countable
    support iterations, generalizing work of 
    Miyamoto, Abraham, and Shelah.
\end{abstract}

\maketitle

%\tableofcontents

\section{Introduction}

In \cite{shelahcountryman}, Shelah proved the existence an uncountable linear order $C$ such that $C \times C$ is a union of countably
many chains in the coordinate-wise partial order,
answering a question
of R. Countryman \cite{countryman}.
Such orders are now known as \emph{Countryman lines} and, as we will see momentarily, are fundamental objects within
the class of linear orders.
Countryman had already observed that such
linear orders are necessarily \emph{Aronszajn}---they have
no uncountable separable suborders and they do not contain a copy of $\omega_1$ or $\omega_1^*$.
Galvin observed that if $C$ is a Countryman line, no uncountable linear order can embed into both $C$ and $C^*$.
Clearly, any uncountable suborder of a Countryman line is Countryman.

Thus, Shelah's result demonstrated that any basis for the Aronszajn lines must contain at least two elements: 
a Countryman line $C$ and its converse $C^*$.
Shelah conjectured \cite[9A]{shelahcountryman} that it was consistent that any Aronszajn line contains a Countryman suborder.
Subsequent work of Abraham and Shelah \cite{AS}, when combined with prior work of Baumgartner \cite{reals_iso}, established
the following reduction:
\begin{thm}[\cite{AS}] \label{SC_equiv}
Assume $\textsf{PFA}$.
The following are equivalent:
\begin{enumerate}
    \item Shelah's conjecture.
    \item If $X \subseteq \Rbb$ has cardinality $\omega_1$ and $C$ is a Countryman line, then any uncountable
    linear order contains an isomorphic copy of $X$, $\omega_1$, $\omega_1^*$, $C$, or $C^*$.
    \item There is an Aronszajn tree $T$ such that $\textsf{CAT}(T)$:
    whenever $K \subseteq T$, there
    is an uncountable antichain $X \subseteq T$ such that $\meet (X)$ is contained in $K$ or $T \setminus K$.
\end{enumerate}
\end{thm}
$\textsf{PFA}$ denotes the {\em Proper Forcing Axiom}, a powerful and versatile axiomatic assumption 
which strengthens Martin's Axiom and the negation
of the Continuum Hypothesis; see \cite{PFA_ICM, comb_dichot}.
In fact, the bounded form of $\textsf{PFA}$ (denoted $\textsf{BPFA}$) introduced by Goldstern and Shelah 
\cite{goldsternshelah} suffices to prove the equivalence in Theorem \ref{SC_equiv}.
$\textsf{BPFA}$ asserts that if a $\Sigma_1$-sentence with parameters in $\mathscr{P}(\omega_1)$ holds in a proper forcing
extension, it is true.
Unlike $\textsf{PFA}$, which has considerable large cardinal strength (see, e.g., \cite{PFA_AD}),
$\textrm{BPFA}$ is equiconsistent with the existence of a reflecting cardinal \cite{goldsternshelah, loc_refl},
a hypothesis which is weaker in consistency strength than the existence of a Mahlo cardinal.

Shelah's conjecture was eventually confirmed by the second author:

\begin{thm}[\cite{linear_basis}] \label{SC}
Assume $\textsf{PFA}$.
Any Aronszajn line contains a Countryman suborder.
\end{thm}

\begin{cor}
    Assume $\textsf{PFA}$ and let $X \subseteq \Rbb$ have cardinality $\omega_1$ and $C$ be any Countryman line.
    Any uncountable linear order contains an isomorphic copy of $X$, $\omega_1$, $\omega_1^*$, $C$, or $C^*$.
\end{cor}

The proof of Theorem \ref{SC} had two unusual related features.
First, it required PFA just to prove that certain posets needed in the argument were proper.
In particular, even though instances of Shelah's conjecture are equivalent to
$\Sigma_1$-sentences with parameters
in $\mathscr{P}(\omega_1)$, the proof of \cite{linear_basis} did not yield that
Bounded Proper Forcing Axiom ($\textsf{BPFA}$) implies Shelah's conjecture. 
Second, the proof relied on consequences of PFA which have large cardinal strength.
Unlike $\textsf{MA}_{\omega_1}$, whose consistency can be established relative to $\textsf{ZFC}$, $\textsf{PFA}$ has considerable large cardinal strength (see \cite{PFA_AD}).
Subsequent work of K\"onig, Larson, Veli\v{c}kovi\'c, and the second author \cite{con_linear_basis} reduced the large cardinal strength of Shelah's conjecture
to something less than a Mahlo cardinal.
Still, this work left open
the possibility that the consistency strength of Shelah's conjecture was stronger than that of $\textsf{BPFA}$,
which is equiconsistent with the existence of a reflecting cardinal \cite{goldsternshelah, loc_refl}.

In the present paper, we remedy these shortcomings.

\begin{thm}
If $A$ is any Aronszajn line and $C$ is any Countryman line, then there is a proper forcing extension
in which $A$ contains an isomorphic copy of $C$ or $C^*$.
\end{thm}

This yields a number of corollaries, which answer questions in the literature.

\begin{cor}
    Let $C$ be any Countryman line and $X \subseteq \Rbb$ have cardinality $\omega_1$.
    If $L$ is any uncountable linear order, there is a proper forcing extension 
    in which $L$ contains an isomorphic copy of $X$, $\omega_1$, $\omega_1^*$, $C$, or $C^*$.
\end{cor}

\begin{cor}
    Assume $\textsf{BPFA}$ and let $X \subseteq \Rbb$ have cardinality $\omega_1$ and $C$ be any Countryman line.
    Any uncountable linear order contains an isomorphic copy of $X$, $\omega_1$, $\omega_1^*$, $C$, or $C^*$.
\end{cor}

\begin{cor}
    The following are equiconsistent:
    \begin{enumerate}
        \item The conjunction of: $\textsf{PFA}$ for posets of cardinality at most $\omega_1$,
        ``Any two $\omega_1$-dense sets of reals are isomorphic,''
        Shelah's conjecture, and $\varphi$.

        \item There is an inaccessible cardinal.
    \end{enumerate}
\end{cor}

In last corollary, $\varphi$ is a combinatorial consequence of $\textsf{PFA}$ which was introduced in
\cite{con_linear_basis}, where it was shown that the conjunction of $\textsf{BPFA}$ and $\varphi$ implies
Shelah's conjecture.

Central to the proof of our main result is
a new forcing iteration preservation lemma for $\omega_1$-trees. 
This lemma says that for any countable support iteration 
$\Seq{ \p_i, \dot \q_i : i \in \beta }$ of proper forcings, where $\beta$ is a limit ordinal, 
and for any $\omega_1$-tree 
$T$ in the ground model, every subtree of $T$ which appears in $V^{\p_\beta}$ 
contains a subtree which lies in $V^{\p_i}$ for some $i \in \beta$. 
We also use this iteration lemma to prove that the statement $\varphi$ is equiconsistent with 
an inaccessible cardinal, which answers a question of \cite{con_linear_basis}. 

\bigskip

\emph{Terminology and notation:} 
A tree $T$ is an \emph{$\omega_1$-tree} if it has height $\omega_1$ 
and countable levels. 
Level $\alpha$ of $T$ is denoted by $T_\alpha$, and $T \restriction \beta = 
\bigcup_{\alpha < \beta} T_\alpha$. 
If $x \in T$ and $\beta$ is less than the height of $x$, we write 
$x \restriction \beta$ for the unique member of $T$ with height $\beta$. 
By a \emph{subtree} of an $\omega_1$-tree we mean a subset which is 
uncountable and downwards closed. 
Whenever $T_0,\ldots,T_{n-1}$ are $\omega_1$-trees, we write $T_0 \otimes \cdots \otimes T_{n-1}$ 
to denote the product tree consisting of $n$-tuples in $T_0 \times \cdots \times T_{n-1}$ 
whose elements all have the same height, ordered componentwise. 
If all of these trees are equal to the same tree $T$, this product 
tree is denoted by $T^{\otimes n}$.

An $\omega_1$-tree $T$ is \emph{coherent} if it is a downwards closed subtree 
of ${}^{<\omega_1} 2$ satisfying that for all $s, t \in T$, 
the set of $\alpha \in \dom(s) \cap \dom(t)$ such that $s(\alpha) \ne s(\beta)$ is finite. 
For any $s, t \in T$, let $s \meet t$ be the greatest lower bound of $s$ and $t$ in $T$. 
For an antichain $X \subseteq T$, let $\meet(X) = \{ \meet(s,t) : s, t \in X \}$. 
For incomparable $s, t \in T$, let $\Delta(s,t)$ be the least $\alpha \in \dom(s) \cap \dom(t)$ 
such that $s(\alpha) \ne t(\alpha)$. 
Note that $\Delta(s,t)$ is the height of the meet $s \meet t$. 
For a set $Z \subseteq T$ and $t \in T$, let 
$\Delta(Z,t)$ be the set of all ordinals of the form 
$\Delta(s,t)$, where $s \in Z$ is incomparable with $t$. 

\section{Saturation of \texorpdfstring{$\omega_1$}{w1}-trees, \texorpdfstring{$\varphi$}{varphi}, and \texorpdfstring{$\psi$}{psi}}

In this section, we review the concept of saturation for $\omega_1$-trees 
and the related statements $\varphi$ and $\psi$. 
We also briefly discuss the ideas of a subtree base and subtree preservation 
which are due to Baumgartner \cite{baumbase}.
 
\begin{defn}[\cite{con_linear_basis}]
	Let $T$ be an $\omega_1$-tree.
	\begin{enumerate}
		\item Subtrees $U$ and $W$ of $T$ are \emph{almost disjoint} if $U \cap W$ 
		is countable.
		\item An \emph{antichain of subtrees of $T$} is a pairwise almost disjoint 
		family of subtrees of $T$.
		\item A family of subtrees $\Fcal$ of $T$ is \emph{predense} if any subtree of $T$ 
		has uncountable intersection with some member of $\Fcal$.
		\item For a given family of subtrees $\Fcal$ of $T$, $\Fcal^{\perp}$ is defined as 
		the set of all subtrees of $T$ which are almost disjoint 
		with every member of $\Fcal$.
	\end{enumerate}
\end{defn}

Note that if $\Fcal$ is a family of subtrees of an $\omega_1$-tree $T$, then 
$\Fcal \cup \Fcal^{\perp}$ is always predense, and $\Fcal$ itself is predense iff 
$\Fcal^{\perp} = \emptyset$.

\begin{defn}[\cite{con_linear_basis}]
	An $\omega_1$-tree $T$ is \emph{saturated} if every antichain of subtrees of $T$ 
	has cardinality less than $\omega_2$.
\end{defn}

\begin{defn}[\cite{con_linear_basis}] \label{define phi}
	Let $T$ be an $\omega_1$-tree and let $\Fcal$ be a family of subtrees of $T$. 
	Define $\varphi(T,\Fcal)$ as the statement that there exists a sequence 
	$\Seq{U_\xi : \xi \in \omega_1}$ of subtrees of $T$ and a club 
	$C \subseteq \omega_1$ such that for all $\alpha \in C$ and for any 
	$t \in T_\alpha$, either:
	\begin{enumerate}
		\item \label{perp_capture} there exists $\xi < \alpha$ such that $t \in U_\xi$ and 
		$U_\xi \in \Fcal^{\perp}$, or
		\item there exists $\beta < \alpha$ such that for all 
		$\nu \in (\beta,\alpha) \cap C$, there is $\xi < \nu$ such that 
		$t \restriction \nu \in U_\xi$ and $U_\xi \in \Fcal$.
	\end{enumerate}
	Define $\psi(T,\Fcal)$ as the same statement except with alternative (\ref{perp_capture}) omitted.
\end{defn}

When $T$ is understood from context, we write $\varphi(\Fcal)$ and $\psi(\Fcal)$ 
to abbreviate $\varphi(T,\Fcal)$ and $\psi(T,\Fcal)$. 
Observe that if $\Fcal$ is predense, then $\varphi(\Fcal)$ and $\psi(\Fcal)$ are equivalent. 
Also note that $\varphi(\Fcal)$ implies $\psi(\Fcal \cup \Fcal^{\perp})$.
It is easily seen that $\varphi(\Fcal)$ and 
$\psi(\Fcal)$ are upwards absolute with respect 
to any $\omega_1$-preserving forcing extension. 

\begin{lem}[{\cite[Lemma 2.2]{con_linear_basis}}] \label{psi implies predense}
	Let $T$ be an $\omega_1$-tree and let $\Fcal$ a family of subtrees of $T$. 
	If $\psi(\Fcal)$ holds, then $\Fcal$ is predense.
\end{lem}

The next two lemmas are implicit in \cite{con_linear_basis}.

\begin{lem} \label{phi implies saturated}
	Suppose that $T$ is an $\omega_1$-tree and for every predense 
	family $\Fcal$ of subtrees of $T$, $\psi(\Fcal)$ holds. 
	Then $T$ is saturated.
\end{lem}

\begin{proof}
	If $T$ is not saturated, then there exists an antichain $\Fcal$ of subtrees of $T$ 
	with cardinality at least $\omega_2$, which we may assume is maximal and 
	hence predense. 
	Let $\Seq{U_\xi : \xi \in \omega_1}$ and $C \subseteq \omega_1$ 
	witness that $\psi(\Fcal)$ holds. 
	Define $\Fcal^* = \{ U_\xi : \xi \in \omega_1 \}$. 
	Note that $\Seq{ U_\xi : \xi \in \omega_1 }$ and $C$ also witness 
	that $\psi(\Fcal^*)$ holds. 
	By Lemma \ref{psi implies predense}, $\Fcal^*$ is predense. 
	Since $\Fcal$ is an antichain, it must be the case that $\Fcal^* = \Fcal$ 
	which is a contradiction.
\end{proof}

\begin{lem} \label{sat phi imply global phi}
	Suppose that $T$ is a saturated $\omega_1$-tree and for any family $\Fcal$ 
	of subtrees of $T$ of size at most $\omega_1$, $\varphi(\Fcal)$ holds. 
	Then for every family $\Fcal$ of subtrees of $T$, $\varphi(\Fcal)$ holds.
\end{lem}

\begin{proof}
	It suffices to show that whenever $\Fcal$ is a family of subtrees of $T$, 
	there exists a subfamily $\Fcal^* \subseteq \Fcal$ of size at most $\omega_1$ 
	such that $\Fcal^{\perp} = (\Fcal^*)^{\perp}$. 
	For then any sequence and club which witness that $\varphi(\Fcal^*)$ holds also 
	witnesses that $\varphi(\Fcal)$ holds. 
	If not, then we can prove that $T$ is not saturated by 
	recursively constructing an almost disjoint 
	sequence $\Seq{ U_\alpha : \alpha \in \omega_2 }$ 
	of subtrees of $T$, each contained in some member of $\Fcal$, as follows. 
	Assume that $\alpha \in \omega_2$ and for all $\beta < \alpha$, 
	$U_\beta$ is defined and contained in some $W_\beta \in \Fcal$. 
	By assumption, $\{ W_\beta : \beta < \alpha \}^{\perp} \ne \Fcal^{\perp}$, so we can 
	find a subtree $U$ of $T$ 
	which is almost disjoint from $W_\beta$ for all $\beta < \alpha$, 
	but for some $W \in \Fcal$, $W \cap U$ is uncountable. 
	Let $U_\alpha = W \cap U$ and $W_\alpha = W$.
\end{proof}

\begin{defn}[\cite{con_linear_basis}]
	\emph{Aronszajn tree saturation} is the statement that every Aronszajn 
	tree is saturated.
\end{defn}

\begin{defn}[\cite{con_linear_basis}]
	Let $\varphi$ denote the statement that $\varphi(T,\Fcal)$ holds for every 
	Aronszajn tree $T$ and any family of subtrees $\Fcal$ of $T$. 
    Define $\psi$ similarly.
\end{defn}

By Lemma \ref{phi implies saturated}, $\varphi$ implies Aronszajn tree saturation.

\begin{rem}
    In \cite{con_linear_basis}, $\varphi$ and $\psi$ were only defined for Aronszajn trees.
    The generalization for $\omega_1$-trees is equivalent, however: 
    if $T$ is an $\omega_1$-tree and $U$ is an Aronszajn tree,
    then $T \otimes U$ is Aronszajn, and if 
    $\varphi(T \otimes U,\Gcal)$ holds for any family of subtrees $\Gcal$ of $T \otimes U$, 
    then $\varphi(T,\Fcal)$ holds for any family of subtrees $\Fcal$ of $T$.
\end{rem}

We now briefly discuss the idea of a subtree base which was introduced in \cite{baumbase}.

\begin{defn}[\cite{baumbase}]
	Let $T$ be an $\omega_1$-tree. 
	A family $\Fcal$ of subtrees of $T$ is a \emph{subtree base} if every subtree 
	of $T$ contains a member of $\Fcal$.
\end{defn}

\begin{defn}[\cite{baumbase}]
	An $\omega_1$-tree has a \emph{small subtree base} if it has a subtree 
	base of size at most $\omega_1$.
\end{defn}

Note that if an $\omega_1$-tree has a small subtree base, then it is saturated. 
Baumgartner and {Todor\v{c}evi\'{c}} \cite{baumbase} proved that if there exists a Kurepa tree, 
then there is a special Aronszajn tree with no small subtree base. 
In fact, if $T$ is a special Aronszajn tree and $U$ is a Kurepa tree, 
then $T \otimes U$ is a special non-saturated Aronszajn tree. 
On the other hand, 
Stejskalová and the first author proved that the existence of a 
non-saturated Aronszajn tree does not imply that there exists a Kurepa tree 
(\cite{jk43, jk46}).

Baumgartner \cite{baumbase} proved that the L\'{e}vy collapse which collapses an inaccessible 
cardinal to become $\omega_2$ forces that every $\omega_1$-tree has a small subtree base. 
The proof is based on the fact that countably closed forcings satisfy the 
property described in the next definition.

\begin{defn}
	For an $\omega_1$-tree $T$, 
	a forcing $\p$ \emph{preserves subtrees of $T$} if $\p$ forces that every 
	subtree of $T$ in $V^\p$ contains a subtree which lies in $V$. 
	A forcing $\p$ \emph{preserves subtrees} if for every $\omega_1$-tree $T$, $\p$ 
	preserves subtrees of $T$.
\end{defn}

The proof of the next lemma is easy.

\begin{lem} \label{preserves subtrees implies predense}
	Let $T$ be an $\omega_1$-tree. 
	Suppose that $\p$ preserves subtrees of $T$. 
	If $\Fcal$ is a predense family of subtrees of $T$, then $\p$ forces 
	that $\Fcal$ is predense.
\end{lem}

\section{Forcing \texorpdfstring{$\varphi$}{varphi} and \texorpdfstring{$\psi$}{psi}}

It was shown in \cite{con_linear_basis} that 
instances of $\varphi$ and $\psi$ can be forced by adding reflecting 
sequences for open stationary set mappings. 
We review this method and prove a new result that the standard forcing 
for adding a reflecting sequence with countable conditions preserves subtrees.

\begin{defn}[\cite{setmappingreflection}] \label{open stationary}
	Let $X$ be an uncountable set and let $M$ be a countable elementary submodel 
	of $H(\theta)$, where $\theta$ is a regular cardinal and $[X]^\omega \in M$. 
	Let $S \subseteq [X]^\omega$.
	\begin{enumerate}
		\item $S$ is \emph{$M$-stationary} if for any club $E \subseteq [X]^\omega$ 
		in $M$, $M \cap S \cap E \ne \emptyset$.
		\item $S$ is \emph{open} if for all $N \in S$ there exists a finite 
		set $x \subseteq N$ such that 
		$$
		\{ P \in [X]^\omega : x \subseteq P \subseteq N \} 
		\subseteq S.
		$$
	\end{enumerate}	
\end{defn}

\begin{defn}[\cite{setmappingreflection}]
	A function $\Sigma$ is an \emph{open stationary set mapping} if there exists 
	an uncountable set $X_\Sigma$ and a regular cardinal $\theta_\Sigma$ such that 
	$\dom(\Sigma)$ is a club subset of $[H(\theta_\Sigma)]^\omega$ consisting of elementary 
	submodels of $H(\theta_\Sigma)$ with $[X_\Sigma]^\omega$ as an member, and for all 
	$M \in \dom(\Sigma)$, $\Sigma(M)$ is open and $M$-stationary.
\end{defn}

\begin{defn}[\cite{setmappingreflection}] \label{mrp forcing}
	Let $\Sigma$ be any open stationary set mapping. 
	Define $\p_{\Sigma}$ as the forcing poset consisting of conditions, ordered 
	by reverse extension, which are 
	$\in$-increasing and continuous sequences $\Seq{ N_i : i \le \alpha }$ 
	of elements of $\dom(\Sigma)$, for some $\alpha \in \omega_1$, 
	such that for all $0 < \beta \le \alpha$ there exists $\xi < \beta$ such that for all 
	$\gamma \in (\xi,\beta)$, $N_\gamma \cap X_\Sigma \in \Sigma(N_\beta)$.
\end{defn}

\begin{thm}[\cite{setmappingreflection}] \label{MRP forcing is proper}
	Let $\Sigma$ be an open stationary set mapping. 
	Then $\p_\Sigma$ is proper and the union of any generic filter on $\p_\Sigma$ 
	is a reflecting sequence for $\Sigma$, by which we mean an 
	$\in$-increasing and continuous sequence 
	$\Seq{ N_\alpha : \alpha \in \omega_1 }$ of elements of $\dom(\Sigma)$
	such that for all non-zero $\alpha \in \omega_1$, there exists $\beta < \alpha$ such that 
	for all $\gamma \in (\beta,\alpha)$, $N_\gamma \cap X_\Sigma \in \Sigma(N_\alpha)$.
\end{thm}

\begin{thm} \label{mrp preserves subtrees}
	Let $\Sigma$ be an open stationary set mapping. 
	Then $\p_\Sigma$ preserves subtrees.
\end{thm}

\begin{proof}
	Let $X = X_\Sigma$ and $\theta = \theta_\Sigma$. 
	Fix an $\omega_1$-tree $T$ and a $\p$-name $\dot U$ for a subtree of $T$. 
	It suffices to prove that for every $p \in \p_{\Sigma}$, there exists $r \le p$ 
	such that the set $\{ x \in T : r \Vdash x \in \dot U \}$ is uncountable. 
	Suppose for a contradiction that $p$ is a counter-example. 
	Fix a large enough regular cardinal $\lambda$ and a countable elementary submodel 
	$M$ of $H(\lambda)$ such that $[X]^\omega$, $\theta$, 
	$\Sigma$, $\p_{\Sigma}$, $H(|\p_\Sigma|^+)$, $p$, $T$, and $\dot U$ 
	are members of $M$. 
	Let $M' = M \cap H(\theta)$. 
	Enumerate $T_{M \cap \omega_1}$ as $\Seq{ t_n : n \in \omega }$ and  
	let $\Seq{ D_n : n \in \omega }$ enumerate all dense open subsets of $\p_\Sigma$ 
	which lie in $M$.

	We construct by recursion a descending sequence of conditions 
	$\Seq{ p_n : n \in \omega }$ in $M \cap \p$ as follows. 
	Let $p_0 = p$. 
	Now let $n \in \omega$ and assume that we have defined $p_n$. 
	Define $E_n$ as the set of all $N \in [X]^\omega$ such that for some 
	countable $N^* \prec H(|\p_{\Sigma}|^+)$ which contains as elements 
	$\Sigma$, $\p_{\Sigma}$, $p_n$, and $D_n$, $N$ is equal to $N^* \cap X$. 
	Then $E_n \in M$ and $E_n$ is a club. 
	Since $\Sigma(M')$ is open and $M'$-stationary, 
	we can fix $N_n \in M' \cap \Sigma(M') \cap E_n$ 
	and $z_n \in [N_n]^{<\omega}$ such that 
	$\{ P \in [X]^\omega : z_n \subseteq P \subseteq N_n \} \subseteq \Sigma(M')$. 
	Let $N_n^* \in M$ witness that $N_n \in E_n$. 
	By the elementarity of $N_n^*$, fix $K_n \in N_n^* \cap \dom(\Sigma)$ such that 
	$z_n \subseteq K_n$. 
	Define $q_n = p_n \cup \{ (\dom(p_n),K_n) \}$, which is in $\p_\Sigma \cap N_n^*$ 
	and extends $p_n$. 
	By assumption and the elementarity of $N_n^*$, 
	there exists $\alpha_n < N_n^* \cap \omega_1$ such that 
	for all $x \in T_{\alpha_n}$, $q_n$ does not force that $x \in \dot U$. 
	In particular, we can find $p_{n+1} \le q_n$ in $N_n^* \cap D_n$ which forces that 
	$t_n \restriction \alpha_n \notin \dot U$. 
	Since $\dot U$ is forced to be downwards closed, 
	$p_{n+1}$ forces that $t_n \notin \dot U$.
	
	This completes the recursion. 
	Define
	$$
	q = \bigcup_n p_n \cup \{ (M \cap \omega_1,M') \}.
	$$
	Note that for all $\gamma \in [\dom(p_0),M \cap \omega_1)$, there exists $n$ such that 
	$$
	z_n \subseteq K_n \cap X \subseteq p_{n+1}(\gamma) \cap X = q(\gamma) \cap X \subseteq N_n,
	$$
	and hence $q(\gamma) \cap X \in \Sigma(M')$. 
	It easily follows that $q \in \p_\Sigma$ and $q \le p_n$ for all $n$. 
	By construction, $q$ forces that $t_n \notin \dot U$ for all $n$. 
	Hence, $q \Vdash \dot U \subseteq \check{T} \restriction (M \cap \omega_1)$, which 
	contradicts that $\dot U$ is forced to be a subtree.
\end{proof}

In \cite{con_linear_basis}, a method for forcing $\varphi(\Fcal)$ using the 
poset of Definition \ref{mrp forcing} was introduced and is based on the following lemma.

\begin{lem}[{\cite[Lemma 3.6]{con_linear_basis}}] \label{Mstationary or perp}
	Let $T$ be an $\omega_1$-tree and let $\Fcal$ be a family of subtrees of $T$. 
	Let $E_0$ be the club subset of $[H(\omega_2)]^\omega$ consisting of all 
	elementary submodels of $H(\omega_2)$ which contain $T$ as a member. 
	Let $\kappa \ge {2^{\omega_1}}^+$ be regular and let 
	$M$ be a countable elementary submodel of $H(\kappa)$ 
	with $T$ and $\Fcal$ members of $M$. 
	Consider any $t \in T_{M \cap \omega_1}$ and suppose that there does not exist 
	a set $A \in M \cap \Fcal^{\perp}$ such that $t \in A$. 
	Let $\Sigma_M$ be the set of all $P \in [H(\omega_2)]^\omega$ such that, if $P \in E_0$, then 
	there exists some $U \in P \cap \Fcal$ such that $t \restriction (P \cap \omega_1) \in U$. 
	Then $\Sigma_M$ is open and $M$-stationary.
\end{lem}

We force $\varphi(\Fcal)$ by decomposing $T$ as a two-dimensional matrix and then force 
witnesses to $\varphi(\Fcal)$ along each column in a forcing iteration of length $\omega$. 
The same construction also forces $\psi(\Fcal)$ in the case that $\Fcal$ is predense.

\begin{prop} \label{forcing phi}
	Let $T$ be an $\omega_1$-tree and let $\Fcal$ be a family of subtrees of $T$. 
	Then there exists a countable support forcing iteration 
	$\Seq{ \p_n, \dot \q_n : n \in \omega }$ satisfying:
	\begin{enumerate}
		\item for each $n \in \omega$, $\p_n$ forces 
		that $\dot \q_n = \p_{\dot \Sigma_n}$ for some open stationary set mapping 
		$\dot \Sigma_n$;
		\item $\p_\omega$ is proper;
		\item $\p_\omega$ forces $\varphi(\Fcal)$;
		\item if $\Fcal$ is predense, then $\p_\omega$ forces $\psi(\Fcal)$.
	\end{enumerate}
\end{prop}

\begin{proof}
	We start by describing the open stationary set mappings used 
	in the factors of the forcing iteration. 
	Suppose that $\{ t_\alpha : \alpha \in \omega_1 \}$ is a set of elements of $T$ 
	such that each $t_\alpha$ has height $\alpha$. 
	Define $\Sigma$ as follows. 
	Let $X_\Sigma = H(\omega_2)$ and $\theta_\Sigma = {2^{\omega_1}}^+$. 
	Define $E_0$ as the club set of all countable elementary submodels of $H(\omega_2)$ 
	which contain $T$ as a member. 
	Let $\dom(\Sigma)$ be the set of all countable elementary submodels of 
	$H(\theta_\Sigma)$ which contain $T$ and $\Fcal$ as members. 
	Consider $M \in \dom(\Sigma)$. 
	If there exists some $A \in M \cap \Fcal^{\perp}$ 
	such that $t_{M \cap \omega_1} \in A$, then let $\Sigma(M)$ 
	be equal to $[H(\omega_2)]^\omega$.
	Otherwise, 	let $\Sigma(M)$ 
	be the set of all $P \in [H(\omega_2)]^\omega$ such that, if $P \in E_0$, then 
	there exists some $U \in P \cap \Fcal$ such that 
	$t \restriction (P \cap \omega_1) \in U$. 
	Then $\Sigma(M)$ is open and $M$-stationary by Lemma \ref{Mstationary or perp}. 
	Now it is straightforward to use a reflecting sequence for $\Sigma$ 
	given by a generic filter on $\p_\Sigma$ to define a sequence 
	$\Seq{ U_\xi : \xi \in \omega_1 }$ which witnesses the 
	definition of $\varphi(\Fcal)$ restricted to elements of 
	$\{ t_\alpha : \alpha \in \omega_1 \}$.
	
	Now working in the ground model, 
	decompose the tree $T$ as $\{ t(\alpha,n) : \alpha \in \omega_1, \ n \in \omega \}$, 
	where each $t(\alpha,n)$ has height $\alpha$. 
	Define a countable support forcing iteration 
	$\Seq{ \p_n, \ \dot \q_n : n \in \omega }$ so that for each $n$, 
	$\p_n$ forces that $\dot \q_n$ equals $\p_{\dot \Sigma_n}$, where $\dot \Sigma_n$ 
	is a $\p_n$-name for the open stationary set mapping described in the previous paragraph 
	for the set $\{ t(\alpha,n) : \alpha \in \omega_1 \}$. 
	In any generic extension by $\p_\omega$, it is easy to combine the $\omega$-many 
	sequences obtained at each stage of the forcing iteration 
	into a single sequence which witnesses that $\varphi(\Fcal)$ holds. 
	In the case that $\Fcal$ is predense, Theorem \ref{mrp preserves subtrees} implies that each 
	$\p_n$ preserves subtrees, and therefore by Lemma \ref{preserves subtrees implies predense} 
    preserves the fact that $\Fcal$ is predense. 
	So each $\p_n$ forces that $\Fcal^{\perp}$ is empty, and consequently the sequence obtained 
	at stage $n+1$ contains only members of $\Fcal$. 
	Hence, the sequence which witnesses $\varphi(\Fcal)$ in $V^{\p_\omega}$ 
	only contains members of $\Fcal$ and therefore witnesses that $\psi(\Fcal)$ holds.
\end{proof}

We note that a different method for forcing $\psi(\Fcal)$ was described 
in \cite{con_linear_basis},  
as Theorem \ref{mrp preserves subtrees} was not known at the time.

\section{An Iteration Lemma for Trees}

In this section, we prove a forcing iteration preservation lemma 
concerning subtrees of $\omega_1$-trees which is an essential 
tool for the main theorems of the article. 
Roughly speaking, we show that for an $\omega_1$-tree $T$ and a countable 
support forcing iteration of proper forcings of limit length, 
any subtree of $T$ which lies in a generic extension by the iteration contains a 
subtree from an intermediate model by an earlier stage of the iteration. 

Let us clarify notation. 
Suppose that $\Seq{ \p_i, \dot \q_i : i \in \beta }$ is a countable support forcing iteration. 
We let $\p_\beta$ denote the countable support limit of this iteration 
in the case that $\beta$ is a limit ordinal, 
and $\p_{\gamma} * \dot \q_\gamma$ in the case that $\beta = \gamma+1$ is a successor ordinal.  
Let $\alpha \in \beta+1$. 
Then $\le_\alpha$ denotes the order on $\p_\alpha$, 
$\dot G_\alpha$ is the canonical $\p_\alpha$-name for the generic filter, and 
$\Vdash_\alpha$ abbreviates the forcing relation $\Vdash_{\p_\alpha}$. 
We also write $\Vdash_\alpha^W$ for the $\p_\alpha$-forcing relation relativized to a model $W$.

Our iteration lemma can be thought of as a natural extension to the theory of 
trees of the following landmark preservation theorem of Shelah. 
This theorem has provided a blueprint for a large number of 
countable support (and revised countable support) 
iteration preservation results across many different areas of set theory. 

\begin{thm}[{\cite[Chapter III, Theorem 3.2]{shelahoriginalbook}}]
	Suppose that $\Seq{ \p_i, \ \dot \q_i : i \in \beta }$ is a 
	countable support forcing iteration of proper forcings. 
	Then $\p_\beta$ is proper.
\end{thm}

One of several possible proofs of Shelah's theorem is based on the following lemma, 
which we will make use of 
(see Lemma 31.17 of \cite{jechbook}).

\begin{lem} \label{standard iteration lemma}
	Let $\Seq{ \p_i, \dot \q_i : i \in \beta }$ 
	be a countable 
	support forcing iteration of proper forcings.  
	Let $\lambda$ be a large enough regular cardinal and let $N$ be a countable 
	elementary submodel of $H(\lambda)$ containing 
	the forcing iteration as a member. 
	Let $\alpha \in N \cap \beta$. 
	Suppose that $\dot p$ is a $\p_{\alpha}$-name for a condition in $\p_\beta$, 
	$q \in \p_\alpha$, and $q \Vdash_{\alpha} \dot p \in N$. 
	Also, assume that $q$ is $(N,\p_\alpha)$-generic and 
	$q \Vdash_{\alpha} \dot p \restriction \alpha \in \dot G_{\alpha}$. 
	Then there exists $q^+ \in \p_\beta$ such that 
	$q^+ \restriction \alpha = q$, $q^+$ is $(N,\p_\beta)$-generic, and 
	$q^+ \Vdash_{\beta} \dot p \in \dot G_\beta$.
\end{lem}

\begin{lem} \label{iteration for trees basic lemma}
	Let $T$ be an $\omega_1$-tree. 
	Suppose that $\Seq{ \p_i, \dot \q_i : i \in \beta }$ is a countable 
	support forcing iteration such that $\p_\beta$ preserves $\omega_1$, 
	$p_0 \in \p_\beta$, 
	and $D$ is a dense open subset of $\p_\beta$. 
	Assume that $p_0$ forces in $\p_\beta$ that $\dot U$ 
	is a downwards closed subset of $\check T$ 
	such that for all $\alpha < \beta$, $\dot U$ contains no subtree 
	which lies in $V^{\p_\alpha}$. 
	Let $\alpha < \beta$. 
	Then $\p_{\alpha}$ forces that, if $p_0 \restriction \alpha \in \dot G_\alpha$, 
	then there exists some $\gamma \in \omega_1$ such that for all $u \in T_\gamma$, 
	there exists some $p_{1,u} \le_\beta p_0$ in $D$ such that 
	$p_{1,u} \restriction \alpha \in \dot G_\alpha$ and 
	$p_{1,u} \Vdash^V_{\beta} u \notin \dot U$.
\end{lem}

\begin{proof}
	Let $r \in \p_{\alpha}$. 
	By extending $r$ further if necessary, we may assume without loss of generality 
	that $r$ decides whether or not $p_0 \restriction \alpha \in \dot G_{\alpha}$, 
	and in the case that it does, that $r \le_\alpha p_0 \restriction \alpha$. 
	If $r$ forces that $p_0 \restriction \alpha \not \in \dot G_{\alpha}$, then we are done. 
	So assume that $r \le_\alpha p_0 \restriction \alpha$. 
	Extend $r \cup p_0 \restriction [\alpha,\beta)$ to some condition $s$ which is in $D$. 
	Then $s \le_\beta p_0$. 
	We claim that $s \restriction \alpha$ forces in $\p_\alpha$ that 
	there exists $\gamma \in \omega_1$ such that for all $u \in T_\gamma$, 
	there exists $p_{1,u} \le_\beta p_0$ in $D$ such that 
	$p_{1,u} \restriction \alpha \in \dot G_\alpha$ and 
	$p_{1,u} \Vdash^V_\beta u \notin \dot U$.
	To prove this claim, 
	let $G_{\alpha}$ be a generic filter for $\p_{\alpha}$ containing $s \restriction \alpha$, 
	and we work in $V[G_{\alpha}]$. 
	Let $\p_\beta / G_{\alpha}$ be the suborder of $\p_\beta$ consisting 
	of all $q \in \p_\beta$ such that $q \restriction \alpha \in G_{\alpha}$. 
	Recall that 
	$\p_\beta$ is forcing equivalent to $\p_{\alpha} * (\p_\beta / \dot G_{\alpha})$. 
	Since $s \restriction \alpha \in G_\alpha$, $s \in \p_\beta / G_\alpha$.

	In $V[G_\alpha]$, let 
	$U_0 = \{ u \in T : s \Vdash_{\p_\beta / G_\alpha} u \in \dot U \}$. 
	Note that $U_0 \in V[G_\alpha]$ and 
	$s$ forces in $\p_\beta / G_\alpha$ that $U_0$ is a downwards closed subset of $\dot U$. 
	Hence, by our assumption on $p_0$ and $\dot U$, $U_0$ is countable. 
	Find $\gamma \in \omega_1$ such that $U_0 \subseteq T \restriction \gamma$. 
	Consider $u \in T_\gamma$. 
	Then $u \notin U_0$, which means that there exists some 
	$s_{u} \le_{\p_\beta / G_\alpha} s$ which forces that $u \notin \dot U$. 
	Then $s_u \le_\beta s$. 
	Since $s_{u} \restriction \alpha \in G_\alpha$, we can find 
	$t_u \le_\alpha s_u \restriction \alpha$ in $G_\alpha$ which forces (in $\p_\alpha$ over $V$) 
	that $s_u$ forces (in $\p_\beta / \dot G_\alpha$ over $V[\dot G_\alpha]$) 
	that $u \notin \dot U$. 
	Let $p_{1,u} = t_u \cup s_u \restriction [\alpha,\beta)$. 
	Then $p_{1,u} \le_\beta p_0$, $p_{1,u} \restriction \alpha = t_u \in G_\alpha$, 
    $p_{1,u} \in D$, and $p_{1,u} \Vdash^V_\beta u \notin \dot U$.
\end{proof}

\begin{lem}[Iteration Lemma for Trees] \label{iteration lemma for trees}
	Let $T$ be an $\omega_1$-tree. 
	Suppose that $\Seq{ \p_i, \ \dot \q_i : i \in \beta }$ is a 
	countable support forcing iteration of proper forcings, where 
	$\beta$ is a limit ordinal. 
	Then $\p_\beta$ forces that for any subtree $U$ of $T$ lying in $V^{\p_\beta}$, 
	there exists some $\alpha < \beta$ and a subtree $W$ in $V^{\p_\alpha}$ 
	such that $W \subseteq U$.
\end{lem}

\begin{proof}
	Let $\dot U$ be a $\p_\beta$-name for a downwards closed subset of $T$. 
	It suffices to show that whenever $p \in \p_\beta$ forces that 
	for all $\alpha < \beta$, $\dot U$ does not contain a subtree which lies in 
	$V^{\p_\alpha}$, then there exists $q \le p$ 
	which forces that $\dot U$ is countable. 
	Without loss of generality, assume that each $\p_i$ is 
	separative, for if not then we can find an equivalent forcing iteration for which 
	this is true.

	Fix a large enough regular cardinal $\lambda$ and let $N \prec H(\lambda)$ 
	be countable which contains the above parameters. 
	Let $\delta = N \cap \omega_1$. 
	Let $\Seq{ D_n : n \in \omega }$ be an enumeration of all of the 
	dense open subsets of $\p_\beta$ which lie in $N$ 
	and let $\Seq{ a_n : n \in \omega }$ 
	be an enumeration of $T_{\delta}$. 
	Fix an increasing sequence $\Seq{ \alpha_n : n \in \omega }$ of ordinals 
	in $N \cap \beta$ cofinal in $\sup(N \cap \beta)$ with $\alpha_0 = 0$.
	
	We define by recursion two sequences $\Seq{ q_i : i \in \omega }$ and 
	$\Seq{ \dot p_i : i \in \omega }$ satisfying that for all $i$:
	\begin{enumerate}
		\item $\dot p_i$ is a $\p_{\alpha_i}$-name for a condition 
		in $\p_\beta$, and $\dot p_0$ is a $\p_0$-name for $p$;
		\item $q_i$ is an $(N,\p_{\alpha_i})$-generic condition 
		and $q_i \Vdash_{\alpha_i} \dot p_i \in N \wedge
		\dot p_i \restriction \alpha_i \in \dot G_{\alpha_i}$;
		\item $q_{i+1} \restriction \alpha_i = q_i$, and 
		$q_{i+1}$ forces in $\p_{\alpha_{i+1}}$ that:
		\begin{enumerate}
			\item $\dot p_{i+1} \le_\beta \dot p_i$;
			\item $\dot p_{i+1} \in D_i$;
			\item $\dot p_{i+1} \Vdash^V_\beta a_i \notin \dot U$.
		\end{enumerate}
	\end{enumerate}

	Begin by letting $q_0$ be the empty condition in $\p_0$ and 
	letting $\dot p_0$ be a $\p_0$-name for $p$. 
	Now let $i \in \omega$ and assume that $\dot p_i$ and $q_i$ are defined as required. 
	Applying Lemma \ref{standard iteration lemma} to $\dot p_i \restriction \alpha_{i+1}$ and $q_i$, 
	find $q_{i+1} \in \p_{\alpha_{i+1}}$ 
	such that $q_{i+1} \restriction \alpha_i = q_i$, $q_{i+1}$ is $(N,\p_{\alpha_{i+1}})$-generic, 
	and $q_{i+1} \Vdash_{\alpha_{i+1}} \dot p_i 
	\restriction \alpha_{i+1} \in \dot G_{\alpha_{i+1}}$.

	To define $\dot p_{i+1}$, consider a generic filter $G_{\alpha_{i+1}}$ 
	which contains $q_{i+1}$. 
	Since $q_{i+1}$ is $(N,\p_{\alpha_{i+1}})$-generic, 
	$N[G_{\alpha_{i+1}}] \cap V = N$. 
	Let $G_{\alpha_i} = G_{\alpha_{i+1}} \cap \alpha_i$. 
	Let $p_i = \dot{p}_i^{G_{\alpha_{i}}}$, which is in $N$. 
	So $p_i \restriction \alpha_{i+1} \in G_{\alpha_{i+1}}$. 
	Also, $p_i \le_\beta p$, so $p_i$ forces (in $\p_\beta$ over $V$) that 
	$\dot U$ is a downwards closed subset of $T$ 
	such that for all $\alpha < \beta$, $\dot U$ has no subtree which lies in $V^{\p_\alpha}$.
	
	Applying Lemma \ref{iteration for trees basic lemma} 
	to $p_i$, $D_i$, and $\alpha_{i+1}$, we conclude that 
	there exists some $\gamma \in \omega_1$ such that for all $u \in T_\gamma$, 
	there is some $p_{1,u} \le_\beta p_i$ in $D_i$ such that 
	$p_{1,u} \restriction \alpha_{i+1} \in G_{\alpha_{i+1}}$ and 
	$p_{1,u} \Vdash^V_\beta u \notin \dot U$. 
	Since $p_i$ and $D_i$ are in $N[G_{\alpha_{i+1}}]$, 
	we may assume that $\gamma < N[G_{\alpha_{i+1}}] \cap \omega_1 = \delta$. 
	Consider $u = a_i \restriction \gamma$, which is in $N[G_{\alpha_{i+1}}]$. 
	By elementarity, we can find 
	$p_{i+1} \in N[G_{\alpha_{i+1}}]$ such that 
	$p_{i+1} \le_\beta p_i$ is in $D_i$, 
	$p_{i+1} \restriction \alpha_{i+1} \in G_{\alpha_{i+1}}$, and 
	$p_{i+1} \Vdash^V_\beta a_i \restriction \gamma \notin \dot U$. 
	As $\dot U$ is forced to be downwards closed, it follows 
	that $p_{i+1} \Vdash^V_{\beta} a_i \notin \dot U$. 
	Finally, since $N[G_{\alpha_{i+1}}] \cap V = N$, $p_{i+1} \in N$.

	Since $G_{\alpha_{i+1}}$ was an arbitrary generic filter 
	on $\p_{\alpha_{i+1}}$ containing $q_{i+1}$, 
	we can fix a $\p_{\alpha_{i+1}}$-name $\dot p_{i+1}$ for a member of $\p_\beta$ 
	which $q_{i+1}$ forces has the properties described above. 
	Now it is routine to check that $q_{i+1}$ and $\dot p_{i+1}$ are as required.

	This completes the construction. 
	Let $q = \bigcup_n q_n$, which is in $\p_{\beta}$. 
	We will argue that $q \le p$ and $q$ forces that $\dot U$ is a subset of $N$, 
	and hence is countable. 
	Let $G_\beta$ be a generic filter 
	which contains $q$ and let $U = \dot U^{G_{\beta}}$. 
	For all $i \in \omega$, let $G_{\alpha_i} = G_\beta \cap \p_{\alpha_i}$ 
	and $p_i = \dot p_i^{G_{\alpha_i}}$. 
	Since $q \restriction \alpha_i$ forces that 
	$\dot p_i \restriction \alpha_i \in \dot G_{\alpha_i}$ 
	and $\p_{\alpha_i}$ is separative, it follows that 
	$q \restriction \alpha_i \le_{\alpha_i} p_i \restriction \alpha_i$. 
	Hence, for all $i < j$, $q \restriction \alpha_j \le_{\alpha_j} 
	p_j \restriction \alpha_j \le_{\alpha_j} 
	p_i \restriction \alpha_j$. 
	Since $\dom(q) \subseteq \sup_n \alpha_n$,  
	it follows that $q \le_\beta p_j$ for all $j$. 
	In particular, $q \le p_0 = p$. 
	Also, $q \in G_\beta$ implies that $p_i \in G_\beta$ for all $i$, and 
	by the choice of $p_{i+1}$, $a_i \notin U$. 
	So $U \subseteq T \restriction \delta \subseteq N$.	
	As $G_\beta$ was arbitrary, $q$ forces all of the above statements. 
	So $q \le p$ and $q \Vdash_{\beta} \dot U \subseteq N$.
\end{proof}

\section{Some Consequences of the Iteration Lemma}

In this section, we derive some consequences of Lemma \ref{iteration lemma for trees}. 
Further applications appear in subsequent sections.

\begin{cor} \label{iteration preserves subtrees}
	Let $T$ be an $\omega_1$-tree. 
	Suppose that $\Seq{ \p_i, \dot \q_i : i \in \beta }$ is a 
	countable support forcing iteration of proper forcing such that for all 
	$i \in \beta$, $\p_i$ forces that $\dot \q_i$ preserves subtrees of $T$. 
	Then $\p_\beta$ preserves subtrees of $T$.
\end{cor}

\begin{proof}
	By induction on $\beta$, using the assumption of the corollary at successor 
	stages and Lemma \ref{iteration lemma for trees} at limit stages.
\end{proof}

Combining Proposition \ref{forcing phi} with Corollary \ref{iteration preserves subtrees}, 
we see that the forcing of Proposition \ref{forcing phi} which adds a $\varphi$-sequence 
preserves subtrees.

\begin{cor} \label{forcing phi 2}
	Let $T$ be an $\omega_1$-tree and let $\Fcal$ be a family of subtrees of $T$. 
	Then there exists a forcing $\p$ which is proper, preserves subtrees, 
	forces $\varphi(\Fcal)$, and if $\Fcal$ is predense, forces $\psi(\Fcal)$.
\end{cor}

\begin{proof}
	Immediate from Theorem \ref{mrp preserves subtrees}, 
	Proposition \ref{forcing phi}, and Corollary \ref{iteration preserves subtrees}.
\end{proof}

\begin{lem} \label{forcing multiple psi}
	Let $\{ (T_i,\Fcal_i) : i \in \delta \}$ be such that for each $i \in \delta$, 
	$T_i$ is an $\omega_1$-tree and $\Fcal_i$ is a predense family of subtrees of $T_i$. 
	Then there exists a proper forcing which preserves subtrees and 
	forces that for all $i \in \delta$, 
	$\psi(T_i,\Fcal_i)$ holds.
\end{lem}

\begin{proof}
	Define a countable support forcing iteration 
	$\Seq{ \p_i, \dot \q_i : i \in \delta }$ as follows. 
	Arrange that each $\p_i$ forces that $\dot \q_i$ preserves subtrees, so 
	by Corollary \ref{iteration preserves subtrees} each $\p_i$ 
	preserves subtrees. 
	Assuming that $\p_i$ is defined, 
	let $\dot \q_i$ be a $\p_i$-name for the 
	forcing of Corollary \ref{forcing phi 2} for $T_i$ and $\Fcal_i$. 
	By Lemma \ref{preserves subtrees implies predense}, 
	$\p_i$ forces that $\Fcal_i$ is predense. 
	Therefore, $\p_i * \dot \q_i$ forces $\psi(T_i,\Fcal_i)$. 
	Since $\psi$ is upwards absolute, $\p_\delta$ is as required.
\end{proof}

We now revisit some classic forcing iteration preservation theorems for trees, and show 
that they can be derived as corollaries of 
Lemma \ref{iteration lemma for trees}. 
Silver \cite{silver} demonstrated the consistency of the non-existence of a Kurepa tree 
from an inaccessible cardinal 
using the fact that countably closed forcings do not add new branches to 
$\omega_1$-trees. 
Shelah proved that for an $\omega_1$-tree $T$, the property of a forcing poset that 
it is proper and does not add new cofinal branches of $T$ is preserved 
under countable support forcing iterations 
(see \cite[Chapter III, Theorem 8.5]{shelahproper}). 
This theorem generalizes Silver's method to a broader 
context by allowing for the construction of models with no Kurepa trees 
using iterated forcing.

Note that any cofinal branch of an $\omega_1$-tree $T$ is a subtree of $T$.

\begin{cor} \label{iteration not adding branches}
	Suppose that $T$ is an $\omega_1$-tree and 
	$\Seq{ \p_i, \dot \q_i : i \in \beta }$ 
	is a countable support forcing iteration of proper forcings 
	such that for all $i \in \beta$, $\p_i$ forces that 
	$\dot \q_i$ does not add new cofinal branches of $T$. 
	Then $\p_\beta$ does not add new cofinal branches of $T$.
\end{cor}

\begin{proof}
	We may assume by induction that for all $i \in \beta$, $\p_i$ forces 
	that any cofinal branch of $T$ in $V^{\p_i}$ lies in $V$. 
	If $\beta$ is a successor ordinal, then we are done by the assumptions of the corollary. 
	Suppose that $\beta$ is a limit ordinal. 
	Then by Lemma \ref{iteration lemma for trees}, 
	any cofinal branch of $T$ in $V^{\p_\beta}$ 
	contains a subtree lying in $V^{\p_i}$ for some $i \in \beta$. 
	But any subtree of a cofinal branch is equal to that branch. 
	By the inductive hypothesis, that cofinal branch is in $V$.
\end{proof}

Abraham and Shelah \cite{AS2} and Miyamoto \cite{miyamoto} independently proved that 
for a given Suslin tree $S$, the 
property of a forcing poset that it is proper and preserves the Suslin property 
of $S$ is preserved by countable 
support forcing iterations. 
This result has a number of significant applications, such as the consistency 
of the forcing axiom $\textsf{PFA}(S)$ for a coherent Suslin tree $S$ (\cite{PFAS}). 
We reprove this result as a special case of Lemma \ref{iteration lemma for trees} 
by making use of the following well-known characterization of 
an $\omega_1$-tree being Suslin.

\begin{lem} \label{suslin}
	Let $S$ be an $\omega_1$-tree. 
	Then $S$ is Suslin iff $S$ contains no cofinal branch and for every subtree 
	$U \subseteq S$, $U$ contains a cone 
	(that is, there exists some $x \in S$ such that 
	the set $\{ y \in S : x \le_S y \}$ is a subset of $U$). 
\end{lem}

\begin{cor}
	Suppose that $S$ is a Suslin tree. 
	Assume that $\Seq{ \p_i, \dot \q_i : i \in \beta }$ 
	is a countable support forcing iteration of proper forcings 
	such that for all $i \in \beta$, $\p_i$ forces that $S$ is Suslin. 
	Then $\p_\beta$ forces that $S$ is Suslin.
\end{cor}

\begin{proof}
	If $\beta$ is a successor ordinal, then the result is immediate from the 
	assumptions of the corollary. 
    Assume that $\beta$ is a limit ordinal. 
	By Corollary \ref{iteration not adding branches}, 
	$\p_\beta$ forces that $S$ does not have a cofinal branch. 
	So it suffices to prove that any subtree of $S$ in $V^{\p_\beta}$ contains a cone. 
	By Lemma \ref{iteration lemma for trees}, $\p_\beta$ forces that 
	any subtree $U$ of $S$ in $V^{\p_\beta}$ contains a subtree $W$ in 
	$V^{\p_i}$ for some $i \in \beta$. 
	Since $S$ is Suslin in $V^{\p_i}$, 
	$W$ contains a cone in $V^{\p_i}$, which by upwards absoluteness 
	implies that $U$ contains a cone in $V^{\p_\beta}$.
\end{proof}

For our final application of Lemma \ref{iteration lemma for trees} 
in this section, 
we show that \textsf{MRP} is consistent with the statement that every 
$\omega_1$-tree has a small subtree base. 
A standard way to obtain a model of \textsf{MRP} is to iterate forcings with countable 
support up to a supercompact cardinal $\kappa$, where each forcing is proper, has 
size less than $\kappa$, and adds a reflecting 
sequence for some open stationary set mapping, 
guided by a Laver function. 
By Theorem \ref{mrp preserves subtrees}, each factor of the iteration preserves 
subtrees, and hence so does the entire iteration by Corollary \ref{iteration preserves subtrees}. 
Let $\Seq{ \p_i : i \in \kappa }$ denote such a forcing iteration. 
By standard arguments, any $\omega_1$-tree $T$ in $V^{\p_\kappa}$ lies in $V^{\p_i}$ 
for some $i \in \kappa$, and since the tail of the iteration in $V^{\p_i}$ 
also preserves subtrees 
by Corollary \ref{iteration preserves subtrees}, 
any subtree of $T$ in $V^{\p_\kappa}$ contains a subtree lying in $V^{\p_i}$. 
As the collection of subtrees of $T$ lying in $V^{\p_i}$ has cardinality $\omega_1$ 
in $V^{\p_\kappa}$, $T$ has a small subtree base in $V^{\p_\kappa}$.

\begin{cor}
	If there is a supercompact cardinal, 
	then there is a proper forcing extension which satisfies the conjunction 
	of \textsf{MRP} and ``all $\omega_1$-trees have a small subtree base.''
\end{cor}

We note that in contrast, $\textsf{MA}_{\omega_1}$ implies that for every 
Aronszajn tree $T$, the minimum size of a subtree base for $T$ is at least $\omega_2$ 
(\cite[Corollary 1.1]{hanazawa}).

\section{The Consistency Strength of \texorpdfstring{$\varphi$}{varphi}}

In this section, we resolve the issue of the consistency strength of $\varphi$. 
In \cite{con_linear_basis}, it was asked whether $\varphi$ implies 
that $\omega_2$ is a Mahlo cardinal in $L$. 
We answer this question negatively as follows.

\begin{thm} \label{inaccessible equiconsistent to phi}
	The following are equiconsistent:
	\begin{enumerate}
		\item There exists an inaccessible cardinal.
		\item \label{inaccessible equiconsistent to phi 2} $\varphi$.
	\end{enumerate}
\end{thm}

\begin{lem} \label{forcing global phi}
	Suppose that $\kappa$ is an inaccessible cardinal,  
	$S \subseteq \kappa$ is stationary, and $\Seq{ s_\alpha : \alpha \in S }$ 
	is a $\Diamond(S)$-sequence. 
	Let $\Seq{ \p_i, \dot \q_i : i \in \kappa }$ be a countable support forcing iteration 
	of proper forcings, where each $\p_i$ has size less than $\kappa$ and $\p_\kappa$ 
	forces that $\check \kappa$ equals $\omega_2$.
	Assume that:
	\begin{enumerate}
		\item \label{force_phi_omega1_stage} $\p_\kappa$ forces that for every $\omega_1$-tree $T$ and 
		for every family $\Fcal$ of subtrees of $T$ of size at most $\omega_1$, 
		$\varphi(T,\Fcal)$ holds;
		\item \label{force_psi_stage}
        for all $\alpha \in S$, if $s_\alpha$ codes (in some canonical way) 
		nice $\p_\alpha$-names for an 
		$\omega_1$-tree $T$ with underlying set $\omega_1$ 
		and a predense family $\Fcal$ of subtrees of $T$, 
		then $\p_{\alpha+1}$ forces $\psi(T,\Fcal)$.
	\end{enumerate}
	Then $\p_\kappa$ forces that for any $\omega_1$-tree $T$ and for any 
	family of subtrees $\Fcal$ of $T$, $\varphi(T,\Fcal)$ holds.
\end{lem}

\begin{proof}
	By Lemma \ref{sat phi imply global phi} and (\ref{force_phi_omega1_stage}), it suffices to show that $\p_\kappa$ 
	forces that every $\omega_1$-tree is saturated.  
	If not, then there are nice $\p_\kappa$-names $\dot T$ and $\dot U_i$ for $i \in \kappa$ 
	and some condition $p \in \p_\kappa$ such that 
	$\p_\kappa$ forces that 
	$\dot T$ is an $\omega_1$-tree with underlying set $\omega_1$, 
	$\{ \dot U_i : i \in \kappa \}$ is a predense family of subtrees of $\dot T$, 
	and if $p$ is in the generic filter, then $\{ \dot U_i : i \in \kappa \}$ 
	is a maximal antichain. 
	By a straightforward argument using the fact that $\p_\kappa$ is $\kappa$-c.c., 
	there exists a club $C \subseteq \kappa$ such that for all $\alpha \in C$:
	\begin{itemize}
		\item $\dot T$ is a nice $\p_\alpha$-name;
		\item for all $\gamma < \alpha$, $\dot U_\gamma$ is a nice $\p_\alpha$-name;
		\item for all $\beta < \alpha$ and for any 
		nice $\p_\beta$-name $\dot W$ for a subtree of $\dot T$, 
		$\p_\alpha$ forces that $\dot W$ has uncountable intersection with some 
		member of $\{ \dot U_\gamma : \gamma < \alpha \}$.
	\end{itemize} 
	Now find $\alpha \in S$ such that $s_\alpha$ codes $\dot T$ and 
	$\{ \dot U_\gamma : \gamma < \alpha \}$.
	
	We claim that $\p_\alpha$ forces that $\{ \dot U_i : i < \alpha \}$ is predense. 
	If not, then there is a condition $q \in \p_\alpha$ and a $\p_\alpha$-name 
	$\dot Y$ for a subtree of $\dot T$ such that for all $\gamma < \alpha$, 
	$q$ forces that $\dot Y \cap \dot U_\gamma$ is countable. 
	Applying Lemma \ref{iteration lemma for trees}, 
	find $r \le q$, $\beta < \alpha$, 
	and a nice $\p_\beta$-name $\dot Z$ for a subtree of $\dot T$ such that $r$ forces 
	that $\dot Z \subseteq \dot Y$. 
	But then $r$ forces that for all $\gamma < \alpha$, 
	$\dot Z \cap \dot U_\gamma$ is countable, 
	which contradicts that $\alpha \in C$. 
	By (\ref{force_psi_stage}), $\p_{\alpha+1}$ forces $\psi(\{ \dot U_\gamma : \gamma < \alpha \})$. 
	By upwards absoluteness, $\p_\kappa$ forces 
	$\psi(\{ \dot U_\gamma : \gamma < \alpha \})$, 
	and hence forces that $\{ \dot U_\gamma : \gamma < \alpha \}$ is predense. 
	This is impossible since $p$ forces that $\{ \dot U_i : i \in \kappa \}$ 
	is an antichain.
\end{proof}

\begin{proof}[Proof of Theorem \ref{inaccessible equiconsistent to phi}]
	The statement $\varphi$ implies Aronszajn tree saturation, which in turn implies 
	the non-existence of a Kurepa tree and hence that $\omega_2$ is inaccessible in $L$. 
	For the converse, assume that $\kappa$ is inaccessible. 
	Then in $L$, $\kappa$ is inaccessible and $\Diamond_\kappa$ holds as witnessed by 
	some diamond sequence $\Seq{ s_\alpha : \alpha \in \kappa }$.  
	Define a countable support forcing iteration $\Seq{ \p_i, \dot \q_i : i \in \kappa }$ 
	of proper forcings by recursion as follows. 
	Arrange that for cofinally many stages of the iteration we force 
	with $\mathrm{Col}(\omega_1,\omega_2)$, thereby ensuring $\p_\kappa$ forces $\check \kappa = \dot \omega_2$. 
	By standard bookkeeping and Proposition \ref{forcing phi}, 
	we can ensure that for any $\omega_1$-tree $T$ and for any family 
	$\Fcal$ of subtrees of $T$ of size at most $\omega_1$ which appear in $V^{\p_\kappa}$,  
	there is some $\alpha \in \kappa$ such that $\p_{\alpha+1}$ forces $\varphi(\Fcal)$. 
	Finally, consider any limit ordinal $\alpha \in \kappa$ such that $s_\alpha$ 
	codes nice $\p_\alpha$-names $\dot T$ and 
	$\{ \dot U_i : i < \alpha \}$ for an $\omega_1$-tree with underlying $\omega_1$ 
	and a predense family of subtrees of $\dot T$. 
	Apply Proposition \ref{forcing phi} to find $\dot \q_\alpha$ 
	such that $\p_\alpha * \dot \q_\alpha$ forces $\psi(\{ \dot U_i : i < \alpha \})$. 
	Now we are done by Lemma \ref{forcing global phi}.
\end{proof}

In the above construction, we can also bookkeep to force along the way 
with any proper forcing of size less than $\kappa$ which appears at some stage less than 
$\kappa$ in the iteration. 
In this way, Lemma \ref{inaccessible equiconsistent to phi}(\ref{inaccessible equiconsistent to phi 2}) 
can be strengthened to include other statements, such as 
the proper forcing axiom for proper forcings of size at most $\omega_1$ or the statement 
that any two $\omega_1$-dense sets of reals are isomorphic. 
In fact, if we strengthen the inaccessibility assumption to a reflecting cardinal, 
this method will produce a model of $\varphi$ together with $\textsf{BPFA}$ 
(\cite{goldsternshelah}). 

\begin{thm}
	The following are equiconsistent:
	\begin{enumerate}
		\item There exists a reflecting cardinal.
		\item $\textsf{BPFA}$ and $\varphi$ hold.
	\end{enumerate}
\end{thm}

\section{The Families \texorpdfstring{$\Rscr_n$}{Rn} and \texorpdfstring{$\Rscr_n^{\perp}$}{Rnperp}}

We now describe families of subtrees of products of a given coherent Aronszajn tree 
which play in a role in the analysis of partial orders involved in forcing an instance of \textsf{CAT}. 
We make use of the 
following basic result about coherent Aronszajn trees.

\begin{lem} \label{fund fact about coherent}
	Let $T$ be a special coherent Aronszajn tree. 
	Assume that $\{ X_\alpha : \alpha \in \omega_1 \}$ is a pairwise disjoint 
	family of members of $T^{\otimes n}$, where $n \in \omega$ is positive. 
	Then there is an uncountable set $B \subseteq \omega_1$ such that for all 
	$\alpha < \beta$ in $B$, every element of $X_\alpha$ is incomparable in $T$ with 
	every element of $X_\beta$, and for all $i < j < n$, 
	$\Delta(X_\alpha(i),X_\beta(i)) = \Delta(X_\alpha(j),X_\beta(j))$.
\end{lem}

See Lemmas 4.2.4 and 4.2.5 of \cite{todorbook} for the proof.

For the remainder of this section, fix a special 
coherent Aronszajn tree $T$ and a set $K \subseteq T$.

\begin{defn}[\cite{linear_basis}]
	Let $n \in \omega$ be positive. 
	For any $X \in T^{\otimes n}$, define $K(X)$ as the set of $\gamma$ less than the 
	height of $X$ such that for all $i < n$, $X(i) \restriction \gamma \in K$.
\end{defn}

\begin{defn}[\cite{con_linear_basis}]
	Let $n \in \omega$ be positive.
	For any uncountable set $Z \subseteq T$, let 
	$R^n_Z$ be the set of all $X \in T^{\otimes n}$ for which 
	there exists some $t \in T$ in the downward closure of $Z$ 
	with the same height as $X$ such that $\Delta(Z,t) \cap K(X) = \emptyset$.
\end{defn}

For any uncountable set $Z \subseteq T$, if $X \in R^n_Z$ as witnessed by $t$, then 
for all $\gamma$ less than the height of $X$, $X \restriction \gamma \in R^n_Z$ 
as witnessed by $t \restriction \gamma$. 
So $R^n_Z$ is downwards closed, and hence is a subtree provided that it is uncountable.

\begin{defn}[\cite{con_linear_basis}]
	Define $\Rscr_n$ as the set of all subtrees of $T^{\otimes n}$ 
	of the form $R^n_Z$ for some uncountable set $Z \subseteq T$.
\end{defn}

\begin{lem} \label{Rnperp implies two in K}
	Let $n \in \omega$ be positive and let $S \in \Rscr_n^{\perp}$. 
	Suppose that $\{ X_\alpha : \alpha \in \omega_1 \}$ 
	is a family of distinct elements of $S$. 
	Then there exist $\alpha < \beta$ such that for all $i < n$, 
	$X_\alpha(i)$ and $X_\beta(i)$ are incomparable in $T$ and 
	$X_\alpha(i) \wedge X_\beta(i) \in K$.
\end{lem}

\begin{proof}
	Suppose for a contradiction that the conclusion of the lemma fails. 
	Applying Lemma \ref{fund fact about coherent}, 
	find an uncountable set $B \subseteq \omega_1$ such that 
	for all $\alpha < \beta$ in $B$ and for all $i < j < n$, 
	$\Delta(X_\alpha(i),X_\beta(i)) = \Delta(X_\alpha(j),X_\beta(j))$. 
	Define $Z = \{ X_\alpha(0) : \alpha \in B \}$. 
	We prove that $\{ X_\xi : \xi \in B \} \subseteq R^n_Z$, which is a contradiction 
	since it then follows that $S \cap R^n_Z$ is uncountable. 
	So let $\xi \in B$. 
	Then $X_\xi(0)$ is in $Z$ and has the same height as $X$, 
	so it suffices to show that 
	$\Delta(Z,X_\xi(0)) \cap K(X_\xi) = \emptyset$. 	
	Consider $\gamma \in \Delta(Z,X_\xi(0))$. 
	Then there is $\alpha \in B$ such that $\gamma = \Delta(X_\alpha(0),X_\xi(0))$. 
	By our assumption for a contradiction, there is $i < n$ such that 
	$X_\alpha(i) \meet X_\xi(i) \not \in K$. 
	By the choice of $B$, $\gamma = \Delta(X_\alpha(i),X_\xi(i))$. 
	So $X_\xi(i) \restriction \gamma = X_\alpha(i) \meet X_\xi(i) \notin K$, and hence 
	$\gamma \notin K(X_\xi)$.
\end{proof}

\begin{lem} \label{nice subtree Rn}
	Let $n \in \omega$ be positive and	 suppose that $U$ is subtree of $T^{\otimes n}$ 
	which is contained in some member of $\Rscr_n$. 
	Then there exists a subtree $S \subseteq U$ such that, letting 
	$Z = \{ Y(0) : Y \in S \}$, for every $X \in S$, 
	$\Delta(Z,X(0)) \cap K(X) = \emptyset$, and in particular, $S \subseteq R^n_Z$.
\end{lem}

\begin{proof}
	Fix an uncountable set $Z_1 \subseteq T$ such that $U \subseteq R^n_{Z_1}$. 
	Choose a sequence $\Seq{ X_\alpha : \alpha \in \omega_1 }$ which lists 
	in increasing heights uncountably many elements of $U$. 
	For each $\alpha \in \omega_1$, fix some $t_\alpha$ in the downward closure of $Z_1$ 
	with the same height as $X_\alpha$ 
	such that $\Delta(Z_1,t_\alpha) \cap K(X_\alpha) = \emptyset$. 
	Applying Lemma \ref{fund fact about coherent} to the set of pairs 
	$\{ (t_\alpha,X_\alpha(0)) : \alpha \in \omega_1 \}$, 
	find an uncountable set $B \subseteq \omega_1$ such that for 
	all $\alpha < \beta$ in $B$, 
	$\Delta(t_\alpha,t_\beta) = \Delta(X_\alpha(0),X_\beta(0))$. 
	Let $S$ be the downward closure in $T^{\otimes n}$ of 
	the set $\{ X_\alpha : \alpha \in B \}$, 
	and define $Z = \{ Y(0) : Y \in S \}$. 	
	Note that $S \subseteq U$.
	
	Consider $X \in S$ and we prove that $\Delta(Z,X(0)) \cap K(X) = \emptyset$. 
	Fix $\alpha \in B$ such that $X \le_{T^{\otimes n}} X_\alpha$. 
	Since $\Delta(Z,X(0)) \subseteq \Delta(Z,X_\alpha(0))$ and 
	$K(X) \subseteq K(X_\alpha)$, it suffices to show that 
	$\Delta(Z,X_\alpha(0)) \cap K(X_\alpha) = \emptyset$. 
	Let $\beta \in B \setminus \{ \alpha \}$ and we claim 
	$\Delta(X_\beta(0),X_\alpha(0)) \notin K(X_\alpha)$. 
	By the choice of $B$, 
	$\Delta(X_\beta(0),X_\alpha(0)) = \Delta(t_\beta,t_\alpha) \in \Delta(Z_1,t_\alpha)$. 
	Since $\Delta(Z_1,t_\alpha) \cap K(X_\alpha) = \emptyset$, it follows that 
	$\Delta(X_\beta(0),X_\alpha(0)) \notin K(X_\alpha)$.
\end{proof}

A combination of Lemmas \ref{iteration lemma for trees} 
and \ref{nice subtree Rn} 
implies the next result.

\begin{lem} \label{previous Rn is predense}
	Suppose that $\Seq{ \p_i, \dot \q_i : i \in \beta }$ is a countable support 
	forcing iteration of proper forcings, where $\beta$ is a limit ordinal. 
	Let $n \in \omega$ be positive. 
	Let $\dot \Fcal$ be a $\p_\beta$-name for the set of elements of  
	$\Rscr_n \cup \Rscr_n^{\perp}$ which lie in $V^{\p_i}$ 
	for some $i \in \beta$. 
	Then $\p_\beta$ forces that $\dot \Fcal$ is predense.
\end{lem}

\begin{proof}
	Consider a generic filter $G$ on $\p_\beta$ and let $\Fcal = \dot \Fcal^G$. 
	In $V[G]$, let $U$ be a subtree of $T^{\otimes n}$. 
	For each $i \in \beta$, define $G_i = G \cap \p_i$. 
	Since $\Rscr_n \cup \Rscr_n^{\perp}$ is predense in $V[G]$, we can find some 
	$W$ in it such that $U \cap W$ is uncountable. 
	First, assume that $W \in \Rscr_n^{\perp}$. 
	Applying Lemma \ref{iteration lemma for trees}, 
	fix $i \in \beta$ 
	and a subtree $S \subseteq U \cap W$ which lies in $V[G_i]$. 
	Since $S \subseteq W$, clearly $S \in \Rscr_n^{\perp}$ in $V[G]$. 
	So $S \in \Fcal$ and $S \cap U = S$ is uncountable.
	
	Secondly, assume that $W \in \Rscr_n$. 
	Applying Lemma \ref{nice subtree Rn}, find in $V[G]$ 
	a subtree $S \subseteq U \cap W$ such that, 
	letting $Z = \{ Y(0) : Y \in S \}$, for every $X \in S$ it is the case that 
	$\Delta(Z,X(0)) \cap K(X) = \emptyset$. 
	Applying Lemma \ref{iteration lemma for trees}, 
	fix $i \in \beta$ and a subtree $S_1 \subseteq S$ which lies in $V[G_i]$. 
	Define $Z^* = \{ X(0) : X \in S_1 \}$, 
	which is an uncountable subset of $Z$ in $V[G_i]$. 
	Let $R = R^n_{Z^*}$, which is in $V[G_i]$. 
	As $S_1 \subseteq S$ and $Z^* \subseteq Z$, it follows that 
    $S_1 \subseteq R$. 
	Now it is easy to check that since $R \in \Rscr_n$ in $V[G_i]$, it is also 
	in $\Rscr_n$ in $V[G]$ (see Lemma \ref{Rn is absolute} below). 
	So $R \in \Fcal$, and $R \cap U$ is uncountable since it contains $S_1$.
\end{proof}

\begin{lem} \label{R1 uncountable}
	Suppose that there exists some $R \in \Rscr_1$ which is uncountable. 
	Then there is an uncountable antichain $X \subseteq T$ such that 
	$\meet(X) \subseteq T \setminus K$.
\end{lem}

\begin{proof}
	Let $Z \subseteq T$ be uncountable such that $R = R^1_Z$. 
	Choose a sequence $\Seq{ (x_\alpha) : \alpha \in \omega_1 }$ with  
	increasing heights consisting of members of $R^1_Z$. 
	For each $\alpha \in \omega_1$, fix some 
	$t_\alpha$ in the downwards closure of $Z$ with the same height as $x_\alpha$ 
	such that $\Delta(Z,t_\alpha) \cap K((x_\alpha)) = \emptyset$. 
	Applying Lemma \ref{fund fact about coherent} to the sequence 
	$\Seq{ (t_\alpha,x_\alpha) : \alpha \in \omega_1 }$, 
	find an uncountable set $B \subseteq \omega_1$ such that 
	for all $\alpha < \beta$ in $B$, 
	$\Delta(t_\alpha,t_\beta) = \Delta(x_\alpha,x_\beta)$. 
	Let $X = \{ x_\alpha : \alpha \in B \}$. 
	We claim that $\meet(X) \subseteq T \setminus K$. 
	So consider $\alpha < \beta$ in $B$. 
	Let $\gamma = \Delta(t_\alpha,t_\beta)$. 
	Since $\gamma \in \Delta(Z,t_\alpha)$ and 
	$\Delta(Z,t_\alpha) \cap K((x_\alpha)) = \emptyset$, 
	$x_\alpha \restriction \gamma \notin K$. 
	By the choice of $B$, $\gamma = \Delta(x_\alpha,x_\beta)$, so 
	$x_\alpha \wedge x_\beta = x_\alpha \restriction \gamma \notin K$.
\end{proof}

\section{Stability for \texorpdfstring{$T$}{T} and \texorpdfstring{$K$}{K}} \label{stability section}

In previous sections we considered a family $\Fcal$ of subtrees of an $\omega_1$-tree, 
together with its orthogonal $\Fcal^{\perp}$, and described forcing $\varphi$ and $\psi$ 
for such a family. 
The situation with families of the form $\Rscr_n$ and $\Rscr_n^{\perp}$ is more 
complicated, since the meaning of these objects can change in forcing 
extensions. 

For the remainder of this section, 
let $T$ be a special coherent Aronszajn tree and let $K \subseteq T$. 
We first note that being a member of $\Rscr_n$ is upwards absolute.

\begin{lem} \label{Rn is absolute}
	Suppose that $\q$ is an $\omega_1$-preserving forcing. 
	Then for any positive $n \in \omega$:
	\begin{enumerate}
	\item \label{Rn_up} For any subtree $R \subseteq T^{\otimes n}$, 
	$R \in \Rscr_n$ implies that $\Vdash_{\q} \check R \in \Rscr_n$. 
	\item \label{Rnperp_down} For any subtree $U \subseteq T^{\otimes n}$, 
	if there exists some $q \in \q$ which forces that $\check U \in \Rscr_n^{\perp}$, then 
	$U \in \Rscr_n^{\perp}$.
	\end{enumerate}
\end{lem}	

\begin{proof}
	(\ref{Rn_up}) It suffices to note that for any uncountable set $Z \subseteq T$, the set 
	$R^n_Z$ is computed the same in both $V$ and $V^\q$. 
	(\ref{Rnperp_down}) follows easily from (\ref{Rn_up}).
\end{proof}

\begin{defn}
	Let $\q$ be an $\omega_1$-preserving forcing. 
	We say that $\q$ is \emph{stable for $T$ and $K$} if for all positive $n \in \omega$, 
	if $U \in \Rscr_n^{\perp}$ then $\Vdash_\q U \in \Rscr_n^{\perp}$. 
	If $q \in \q$ and there exists some $U \in \Rscr_n^{\perp}$ such that 
	$q \Vdash_{\q} \check U \notin \Rscr_n^{\perp}$, 
	then we say that $q$ \emph{destabilizes $T$ and $K$}.
\end{defn}

\begin{lem} \label{subtree preservation implies stable}
	Suppose that $\q$ is an $\omega_1$-preserving forcing which preserves subtrees. 
	Then $\q$ is stable for $T$ and $K$.
\end{lem}

\begin{proof}
	Suppose for a contradiction that $n \in \omega$ is positive, 
	$U \in \Rscr_n^{\perp}$, $G$ is a generic filter on $\q$, and $V[G]$ satisfies that  
	$U \notin \Rscr_n^{\perp}$. 
	Then in $V[G]$, there is some $R \in \Rscr_n$ such that 
	the subtree $U \cap R$ is uncountable. 
	Applying Lemma \ref{nice subtree Rn}, find in $V[G]$ 
	a subtree $S \subseteq U \cap R$ such that, 
	letting $Z = \{ Y(0) : Y \in S \}$, for every $X \in S$ it is the case that 
	$\Delta(Z,X(0)) \cap K(X) = \emptyset$. 
	By the fact that $\q$ preserves subtrees, 
	fix a subtree $W \subseteq S$ which lies in $V$. 
	Define $Z^* = \{ X(0) : X \in W \}$, which is an uncountable subset of $Z$ in $V$. 
	The fact that $W \subseteq S$ and $Z^* \subseteq Z$ easily 
	imply that $W \subseteq R^n_{Z^*}$. 
	So in $V$, $U \cap R^n_{Z^*}$ contains the uncountable set $W$, and hence is uncountable. 
	So $U \notin \Rscr_n^{\perp}$ in $V$, which is a contradiction.
\end{proof}

\begin{defn}
	Let $\varphi(T,K)$ be the statement that for all positive $n \in \omega$, 
	$\varphi(\Rscr_n)$ holds.  
\end{defn}

\begin{lem} \label{forcing phiTK}
	There exists a proper forcing poset which preserves subtrees and forces $\varphi(T,K)$. 
\end{lem}

\begin{proof}
	Define a countable support forcing iteration 
	$\Seq{ \p_n, \dot \q_n : n \in \omega }$ 
	so that for all $n \in \omega$, 
	$\p_n$ forces that $\dot \q_{n}$ is the subtree preserving proper 
	forcing of Corollary \ref{forcing phi 2} 
	which forces $\varphi(\Rscr_{n+1})$. 
	By Corollary \ref{iteration preserves subtrees}, 
	$\p_\omega$ preserves subtrees. 
	To see that $\p_\omega$ forces that 
	$\varphi(T,K)$ holds, consider a generic filter $G$ on $\p_\omega$. 
	Let $n \in \omega$ be positive and let $G_{n} = G \cap \p_{n}$. 
	In $V[G_{n}]$, fix a sequence 
	$\Seq{ U_\xi : \xi \in \omega_1 }$ and a club $C \subseteq \omega_1$ 
	which witness that $\varphi(\Rscr_n)$ holds. 
	By Lemma \ref{Rn is absolute}(\ref{Rn_up}), 
	the same sequence and club witness that $\varphi(\Rscr_n)$ holds in $V[G]$ 
	provided that the elements of this sequence which are in 
	$\Rscr_n^{\perp}$ in $V[G_{n}]$ are still in $\Rscr_n^{\perp}$ in $V[G]$. 
	In other words, $\varphi(\Rscr_n)$ holds in $V[G]$ provided that the quotient forcing 
	$\p_{\omega} / G_{n}$ is stable for $T$ and $K$ in $V[G_{n}]$. 
	But this fact holds by Corollary \ref{iteration preserves subtrees} 
	applied to the tail of the iteration together with 
	Lemma \ref{subtree preservation implies stable}.
\end{proof}

In the next section, we prove that it is consistent that every proper 
forcing is stable for $T$ and $K$. 
We also need to know that this property can be preserved by certain forcings.

\begin{lem} \label{preserving stable}
	Suppose that every proper forcing is stable for $T$ and $K$. 
	Let $\Fcal_n = \Rscr_n \cup \Rscr_n^{\perp}$ for each positive $n \in \omega$.  
	Assume that $\p$ is a proper forcing which forces that $\psi(\Fcal_n)$ holds for 
	all positive $n \in \omega$. 
	Then $\p$ forces that every proper forcing is stable for $T$ and $K$.
\end{lem}

\begin{proof}
	Let $\dot \q$ be a $\p$-name for a proper forcing. 
	Let $G$ be a generic filter on $\p$ and let $H$ be a $V[G]$-generic filter on 
	$\q = \dot \q^G$. 
	Suppose for a contradiction that for some positive $n \in \omega$ and 
	some $S \in \Rscr_n^{\perp}$ in $V[G]$, 
	$S$ is not in $\Rscr_n^{\perp}$ in $V[G][H]$. 
	Fix $R \in \Rscr_n$ in $V[G][H]$ such that $S \cap R$ is uncountable. 	
	By upwards absoluteness, $\psi(\Fcal_n)$ holds in $V[G][H]$, and hence 
	$\Fcal_n$ is predense in $V[G][H]$. 
	So we can fix $U \in \Fcal_n$ such that $U \cap S \cap R$ is uncountable. 
	We claim that $U \in \Rscr_n^{\perp}$ in $V$. 
	If not, then since $\Fcal_n = \Rscr_n \cup \Rscr_n^{\perp}$, 
	$U \in \Rscr_n$ in $V$. 
	By Lemma \ref{Rn is absolute}, $U \in \Rscr_n$ in $V[G]$. 
	But $U \cap S$ is uncountable, which contradicts that 
	$S \in \Rscr_n^{\perp}$ in $V[G]$.
	Now in $V[G][H]$, $U \cap S \cap R \subseteq U \cap R$, and hence 
	$U \cap R$ is uncountable. 
	Therefore, $U \notin \Rscr_n^{\perp}$ in $V[G][H]$. 
	So in $V$, $\p * \dot \q$ is a proper forcing which is not stable 
	for $T$ and $K$, which contradicts our assumptions.
\end{proof}

Starting with a model in which all proper forcings are stable for $T$ and $K$ 
and letting $\Fcal_n = \Rscr_n \cup \Rscr_n^{\perp}$ for all positive $n \in \omega$, 
we can use the forcing of Lemma \ref{forcing multiple psi} to 
obtain a generic extension in which $\psi(\Fcal_n)$ holds for all $n$. 
Since $\psi$ is upwards absolute, any further proper forcing extension also satisfies the same. 
Applying Lemma \ref{preserving stable}, we can in this way 
make the property ``all proper forcings are stable for $T$ and $K$'' 
indestructible with respect to proper forcings.

\section{All Proper Forcings Can Be Stable}

In this section, we prove that it is consistent that all proper forcings are stable.

\begin{thm} \label{stabilize}
    Suppose that $T$ is a special coherent Aronszajn tree and $K \subseteq T$.
    There is a proper poset which forces ``every proper forcing is 
    stable for $T$ and $K$.''
\end{thm}

\begin{proof}
    We begin by recursively constructing $\dot \q_\alpha$ for each ordinal $\alpha$ 
    maintaining that for each $\alpha$, $\Seq{\dot \q_\xi : \xi \in \alpha}$
    defines a countable support forcing iteration 
    $\Seq{\p_\xi,\dot \q_\xi : \xi \in \alpha}$ 
    so that $\dot \q_\alpha$ is a $\p_\alpha$-name for a proper forcing.
	We simultaneously recursively define ordinals $\gamma_\alpha$ and 
    $\p_\alpha$-names $\dot \Fcal_\alpha^n$ for each $n \in \omega$. 
    Since the construction describes a proper class, 
	we will be careful to avoid using global choice to ensure that the definition 
	works within the framework of \textsf{ZFC}.

    If $\Seq{\dot \q_\xi : \xi \in \alpha}$ has been defined and
    $\p_\alpha$ is the countable support iteration, define
    $\dot \Fcal_\alpha^n$ for each $n$ as follows.
    $\dot \Fcal_\alpha^n$ consists of all $(\dot S,p)$ such that:
    \begin{itemize}
        \item $p \in \p_\alpha$, 
        \item $\dot S$ is a nice $\p_\alpha$-name for a subset of $T^{\otimes n}$ 
        and if $\alpha$ is a limit ordinal, then moreover
        $\dot S$ is a nice $\p_\beta$-name for some $\beta < \alpha$,
        and        
        \item $p \forces_{\alpha} \dot S \in \dot \Rscr_n \cup \dot \Rscr_n^\perp$.
    \end{itemize}
    We note that if $\p_\alpha$ is a complete suborder of a larger
    forcing $\p'$, we will also regard $\dot \Fcal_\alpha^n$ as a
    $\p'$-name; with this convention $\dot \Fcal_\alpha^n$
    names the same set in both the intermediate extension by $\p_\alpha$ 
    and in the full generic extension. 
	Define $\gamma_\alpha$ to be the least limit ordinal
    such that:
    \begin{itemize}
        \item there is a $\p_\alpha$-name $\dot \q$ of rank less than $\gamma_\alpha$ 
        for a proper poset which forces $\forall n \, \psi({\dot \Fcal}_\alpha^n)$  
        and 
        \item 
    	if there is a $\p_\alpha$-name $\dot \q$ for a proper forcing
	    such that there is
    	$p \in \p_\alpha$, $n \in \omega$, and a $\p_\alpha$-name $\dot S$ 
    	for a member of $\dot \Rscr_n^\perp$ such that 
    	$$
    	p \forces_{\p_\alpha * \dot \q} \dot S \not \in \dot \Rscr_n^\perp 
    	\land
    	\forall m \, \psi({\dot \Fcal}_\alpha^m)
    	$$
    	then there is such a $\dot \q$ of rank less than $\gamma_\alpha$.
    \end{itemize}
    Roughly speaking, we define $\dot \q_\alpha$ to be the $\p_\alpha$-name 
    for the lottery sum of all proper partial orders
    having rank less than $\gamma_\alpha$ and which force 
    $\forall n \, \psi(\dot \Fcal_\alpha^n)$. 
    Notice that by Lemma \ref{previous Rn is predense}, each 
    $\dot \Fcal_\alpha^n$ is forced by $\p_\alpha$ to be predense and
    hence by Lemma \ref{forcing multiple psi}, this lottery sum is 
    nontrivial.

	More formally, define $\dot \q_\alpha$ to consist of $\check \one$ together with 
     all $p \in \p_\alpha$ and all $(\dot x,p)$ such that:
     \begin{itemize}
         \item $\dot x$ is a $\p_\alpha$-name of rank less than $\gamma_\alpha$ which
         is forced to be an ordered pair whose first coordinate is a poset and 
         whose second element is an element of this poset;
         \item if $\dot \q$ is a $\p_\alpha$-name for the first coordinate of $\dot x$ and
         $\dot q$ is a $\p_\alpha$-name for the second coordinate,
         then $p \forces_{\alpha} \dot q \in \dot \q$;
     \end{itemize} 
     Define $\leq_{\dot \q_\alpha}$ to consist of all $(\dot z,p)$ such that:
     \begin{itemize}
         \item $p \in \p_\alpha$,
         \item $\dot z$ is a $\p_\alpha$-name which has rank less than $\gamma_\alpha$;
         \item $p$ forces in $\p_\alpha$ that $\dot z$ is an ordered pair 
         of elements of $\dot \q_\alpha$;
         \item letting $\dot x$ and $\dot y$ be $\p_\alpha$-names 
         for the first and second coordinates of $\dot z$, respectively,
         either:
         \begin{itemize}
             \item $p$ forces the second coordinate of $\dot x$ is $\check \one$ or
             \item for some $\dot \q$, $\leq_{\dot \q}$, $\dot q_0$, and $\dot q_1$,
             $$((\dot \q,\leq_{\dot \q},\dot q_0),p),((\dot \q,\leq_{\dot \q},\dot q_1),p) 
             \in \dot \q_{\alpha},$$
             $$p \forces_{\alpha} \dot q_1 \leq \dot q_0$$
         \end{itemize} 
     \end{itemize}

    This completes the recursive construction.

    Notice that since properness is preserved by countable support iterations, 
    each $\p_\alpha$ is proper.
    The class of ordinals $\delta$ such that if $\alpha < \delta$, 
    then $\gamma_\alpha < \delta$ is a closed unbounded class.
    Let $\theta$ be the least element of this class of uncountable cofinality.
    We first show that $\p_\theta$ has the property that if $\dot \q$ is a $\p_\theta$-name 
    for a proper partial order
    which forces $\forall n \, \psi({\dot \Fcal}_{\theta}^n)$, 
    then $\p_\theta$ forces that $\dot \q$ is stable for $\check T$ and $\check K$. 
    Then we will prove that $\p_{\theta+1}$ satisfies the conclusion of the lemma.
    
    Toward our first goal, suppose for contradiction that there is a $p \in \p_\theta$ and a 
    $\p_\theta$-name $\dot \q$ for a proper poset which forces
    $\forall n \, \psi({\dot \Fcal}_\theta^n)$ such that $p$ forces that $\dot \q$ 
    is not stable for $\check T$ and $\check K$.
    By extending $p$ if necessary, 
    we may assume that for some positive $n$ and $\dot q$ with 
    $p * \dot q \in \p_\theta * \dot \q$,
    there is a $\p_\theta$-name $\dot S$ such that
    $$
    p \forces_{\theta} \dot S \in \dot \Rscr_n^\perp
    $$
    $$
    p * \dot q \forces_{\p_\theta * \dot \q} \dot S \not \in \dot \Rscr_n^\perp.
    $$
    Let $\dot U$ be a $\p_\theta * \dot \q$-name
    such that $p * \dot q$ forces that 
    $\dot U$ is a subtree of $\dot S$ which is contained in an 
    element of $\dot \Rscr_n$.
    Since $\psi({\dot \Fcal}_\theta^n)$ implies that ${\dot \Fcal}_\theta^n$ is predense,
    we can find $p' * \dot q' \leq p * \dot q$ in $\p_\theta * \dot \q$, an $\alpha < \theta$, 
    and a nice $\p_\alpha$-name $\dot S'$ such that 
    $$
    p' \forces_{\theta} \dot S' \in \dot \Rscr_n \cup \dot \Rscr_n^\perp
    $$
    $$p' * \dot q' \forces_{\p_\theta * \dot \q} \dot S' \cap \dot U \textrm{ is uncountable}.$$
%    for an element of $\dot \Rscr_n \cup \dot \Rscr_n^\perp$ such that 
%    $$p' \forces_{\p_\theta} \dot S' \subseteq \dot S$$
    Since $\mathrm{cof}(\theta) > \omega$,
    by increasing $\alpha$ if necessary we may assume that it contains the support of $p'$.
    
    Notice that since $p' * \dot q' \forces_{\p_\theta * \dot \q} \dot S' \cap \dot S$
    is uncountable, it follows that
    $$p' \forces_{\theta} \dot S' \in \dot \Rscr_n^\perp.$$
    Also we have that
    $$p' * \dot q' \forces_{\p_\theta * \dot \q} \dot S' \not \in \dot \Rscr_n^\perp.$$
    Since membership to $\Rscr_n^\perp$ is downward absolute by Lemma \ref{Rn is absolute},
    $$p' \forces_{\alpha} \dot S' \in \dot \Rscr_n^\perp.$$
    Let $\dot G_\alpha$ be the canonical $\p_\alpha$-name for the generic filter and let 
    $\dot \q'$ be the $\p_\alpha |_{p'}$-name for the poset 
    $(\p_\theta |_{p'}/\dot G_\alpha) * (\dot \q |_{\dot q'})$.
    Then $(\p_\alpha |_{p'}) * \dot \q'$ is forcing equivalent to 
    $(\p_\theta |_{p'}) * (\dot \q |_{\dot q'})$ and hence
    $$p' \forces_{\p_\alpha * \dot \q'} \dot S' \not \in \dot \Rscr_n^\perp.$$
    Also, since $\dot \q_\alpha$ is chosen so that it will force 
    $\forall m \, \psi(\dot \Fcal_\alpha^m)$,
    it is forced that $(\p_\theta|_{p'}/\dot G_\alpha) * \dot \q$ is a proper poset 
    which forces $\forall m \, \psi(\dot \Fcal_\alpha^m)$.
    By definition of $\gamma_\alpha$, there is a $\p_\alpha$-name $\dot \q''$ of rank less 
    than $\gamma_\alpha$ such that:
    $$
    p' \forces_{\alpha} \dot \q'' \
    \textrm{is a proper poset which forces} \ \dot S' \not \in \dot \Rscr_n^\perp 
    \ \textrm{and} \ \forall m \, \psi(\dot \Fcal_\alpha^m).
    $$
    Now extend $p'$ to the condition $p''$ in $\p_{\alpha+1}$ 
    whose $(\alpha+1)\St$-coordinate is
    $$((\dot \q'',\leq_{\dot \q''}),\one_{\dot \q''}).$$
    Clearly, $$p'' \forces_{\alpha+1} \dot S' \not \in \dot \Rscr_n^\perp,$$
    and therefore
    $$p'' \forces_{\theta} \dot S' \not \in \dot \Rscr_n^\perp,$$
    contrary to $p'' \leq p'$ and $p' \forces_{\theta} \dot S' \in \dot \Rscr_n^\perp$. 
    This contradiction completes the proof that 
    $\p_\theta$ forces that any proper forcing which forces 
    $\forall n \, \psi(\dot \Fcal_\theta^n)$ is stable for $\check T$ and $\check K$.
    
    Finally, we will verify that $\p_{\theta+1}$ forces that if $\dot \q$ is a proper poset,
    then $\dot \q$ is stable for $\check T$ and $\check K$. 
    First observe that by upwards absoluteness, 
    $\p_\theta$ forces that $\dot \q_\theta * \dot \q$ is a proper poset 
    which forces $\forall n \, \psi({\dot \Fcal}_\theta^n)$. 
    Thus, $\p_\theta$ forces that $\dot \q_\theta * \dot \q$ is stable for 
    $\check T$ and $\check K$. 
   
    Now suppose for a contradiction
    that $n \in \omega$, 
    $\dot S$ is a $\p_{\theta+1}$-name, and $p * \dot q \in \p_{\theta+1} * \dot \q$ 
    is such that
    $$
    p \forces_{\theta+1} \dot S \in \dot \Rscr_n^\perp
    $$
    $$
    p * \dot q \forces_{\p_{\theta+1} * \dot \q} \dot S \not \in \dot \Rscr_n^\perp
    $$
    Let $\dot U$ be a $\p_{\theta+1} * \dot \q$-name for a subtree 
    of $\dot S$ which $p * \dot q$ forces is contained in an
    element of $\dot \Rscr_n$.
    Since $\p_{\theta+1} * \dot \q$ forces $\dot \Fcal_\theta^n$ is predense,
    there is a $\p_\theta$-name $\dot S'$ for an element of ${\dot \Fcal}_\theta^n$
    and a $p' * \dot q' \leq p * \dot q$ such that
    $$
    p' * \dot q' \forces_{\p_{\theta+1} * \dot \q} \dot S' \cap \dot U \textrm{ is uncountable}.
    $$
    In particular,
    $$
    p' * \dot q' \forces_{\p_{\theta+1} * \dot \q} \dot S' \cap \dot S \textrm{ is uncountable}.
    $$
    Since $p' \forces_{\theta+1} \dot S \in \Rscr_n^\perp$,
    it must be that $p' \restriction \theta \forces_{\theta}
    \dot S' \in \Rscr_n^\perp$.
    But now $\dot U$ witnesses that
    $$p' * \dot q' 
    \forces_{\p_\theta * \dot \q_\theta * \dot \q} \dot S' \not \in \dot \Rscr_n^\perp.$$
    This contradicts our earlier observation that
    $\p_\theta$ forces that $\dot \q_\theta * \dot \q$ is stable for $\check T$ and $\check K$.
    \end{proof}

\section{Rejection and the Key Lemma}

The remaining tools which we need for proving the main theorems of the article are the 
forcing posets of \cite{linear_basis} related to forcing an instance of \textsf{CAT}. 
These forcings involve the notion of rejection and the associated Key Lemma. 
We now review these ideas and their connection with the statement $\varphi$. 
For the remainder of this section and the next, 
fix a special Aronszajn tree $T$ which is coherent, uniform, and binary, 
and fix a set $K \subseteq T$.

\begin{defn} \label{E0 and Ecal}
	Define $E_0$ as the club subset of $[H(\omega_2)]^\omega$ consisting of all 
	countable elementary submodels of 
	$H(\omega_2)$ which contain $T$ and $K$ as elements. 
	Let $\Ecal$ be the family of all club subsets of $[H(\omega_2)]^\omega$ 
	which are subsets of $E_0$.
\end{defn}

\begin{defn}[\cite{linear_basis}]
	Let $n \in \omega$ be positive and let $X \in T^{\otimes n}$. 
	For any $P \in E_0$, we say that $P$ \emph{rejects} $X$ 
	if $X$ has height at least $P \cap \omega_1$ 
	and there exists an uncountable set 
	$Z \subseteq T$ in $P$ and some $t \in T_{P \cap \omega_1}$ 
	in the downward closure of $Z$ such that $\Delta(Z,t) \cap K(X) = \emptyset$.
\end{defn}

In the above, since $\Delta(Z,t) \subseteq P \cap \omega_1$, 
easily $P$ rejects $X$ iff $P$ rejects $X \restriction (P \cap \omega_1)$. 
Note that $P$ does not reject the empty tuple.

The next lemma connects the idea of rejection 
with the families $\Rscr_n$, for positive $n \in \omega$. 
The proof is straightforward.

\begin{lem}[{\cite[Lemma 2.6]{con_linear_basis}}] \label{reject iff R}
	Let $n \in \omega$ be positive. 
	Let $P \in E_0$ and let $X \in T^{\otimes n}$ have height at least $P \cap \omega_1$. 
	Then $P$ rejects $X$ iff there exists some $R \in P \cap \Rscr_n$ such that 
	$X \restriction (P \cap \omega_1) \in R$.
\end{lem}

\begin{lem} \label{stationarily many DNR}
	Let $n \in \omega$ be positive. 
	Suppose that $S$ is a subtree of $T^{\otimes n}$ and there exist 
	stationarily many $P \in E_0$ such that $P$ does not reject any member of $S$. 
	Then $S \in \Rscr_n^{\perp}$.
\end{lem}

\begin{proof}
	Suppose for a contradiction that $S \notin \Rscr_n^{\perp}$. 
	Then there exists $R \in \Rscr_n$ such that $S \cap R$ is uncountable. 
	By the stationarity assumption, find $P \in E_0$ such that $R \in P$ 
	and $P$ does not reject any member of $S$. 
	As $S \cap R$ is uncountable and downwards closed, 
	we can fix some $X \in S \cap R$ with height $P \cap \omega_1$. 
	By Lemma \ref{reject iff R}, $P$ rejects $X$, which contradicts the fact that $X \in S$.
\end{proof}

For any finite set $X \subseteq T$ and any $\gamma \in \omega_1$, let $X \restriction \gamma$ denote the tuple which enumerates the set 
$\{ x \restriction \gamma : x \in X, \ \gamma \le \text{ht}_T(x) \}$ 
in lexicographically increasing order.

\begin{defn}
	Let $P \in E_0$ and let $X \subseteq T$ be finite. 
	We say that $P$ \emph{rejects} $X$ if $X \restriction (P \cap \omega_1)$ 
	is a non-empty tuple which $P$ rejects.
\end{defn}

\begin{defn}[\cite{linear_basis}] \label{key lemma}
	The \emph{Key Lemma} is the statement that for any countable elementary submodel 
	$M$ of $H({2^{2^{\omega_1}}}^+)$ which contains $T$ and $K$ as elements and for 
	any finite set $X \subseteq T$, there exists $E \in M \cap \Ecal$ such that 
	either every element of $M \cap E$ rejects $X$, or no element of $M \cap E$ 
	rejects $X$.
\end{defn}

In \cite{linear_basis}, it is shown that \textsf{MRP} implies the Key Lemma. 

\begin{lem} \label{phiTK implies key lemma}
	Assume that $\varphi(T,K)$ holds. 
	Then the Key Lemma holds.
\end{lem}

For a proof, see Lemma 2.7 of \cite{con_linear_basis}.

\section{Instances of \textsf{CAT}}

Our final step before proving the main theorems of the article is to describe 
the forcings of \cite{linear_basis} which are related to forcing an instance of \textsf{CAT}, and 
how the analysis of these forcings in \cite{linear_basis} is  
connected with the idea of stability introduced in Section \ref{stability section}.

Let $\chi$ denote $({2^{2^{\omega_1}}})^+$ for the rest of the section. 
The following definitions all appear in \cite{linear_basis}.

\begin{defn}
	Define $E_1$ as the set of all $N \in [H({2^{\omega_1}}^+)]^\omega$ 
	such that for some countable $M \prec H(\chi)$ with $T$ and $K$ in $M$, 
	$N = M \cap H({2^{\omega_1}}^+)$. 
\end{defn}

Note that $E_1$ is a club.

\begin{defn}
	Let $\Hcal (K)$ denote the set of all finite antichains $X \subseteq T$ 
	such that $\meet(X) \subseteq K$.
\end{defn}

\begin{defn}
	Define $\partial (K)$ as the forcing poset consisting of all pairs 
	$p = (X_p,\Ncal_p)$ such that:
	\begin{enumerate}
		\item $\Ncal_p$ is a finite $\in$-chain of members of $E_1$;
		\item $X_p \subseteq T$ is a finite antichain;
		\item for all $N \in \Ncal_p$, there exists some $E \in N \cap \Ecal$ 
		such that $X_p$ is not rejected by any member of $N \cap E$.
	\end{enumerate}
	Let $q \le p$ if $\Ncal_p \subseteq \Ncal_q$ and $X_p \subseteq X_q$.
\end{defn}

\begin{defn}
	Let $\partial \Hcal (K)$ be the suborder of $\partial (K)$ consisting 
	of conditions $p$ such that $X_p \in \Hcal (K)$.
\end{defn}

Note that $\partial (K)$ and $\partial \Hcal (K)$ are both in $H(\chi)$.

\begin{defn}
	The forcing $\partial \Hcal (K)$ is \emph{canonically proper} if whenever 
	$M$ is a countable elementary submodel of 
	$H(\chi)$ with $T$, $K$, and $\partial \Hcal (K)$ members of $M$, 
	if $p \in \partial (K)$ and $M \cap H({2^{\omega_1}}^+)$ is in $\Ncal_p$, 
	then $p$ is $(M,\partial (K))$-generic. 
	A similar definition is made for $\partial (K)$ being canonically proper.
\end{defn}

\begin{lem}[{\cite[Lemma 5.20]{linear_basis}}] \label{partialHK done}
	Assume that $\partial \Hcal (K)$ is canonically proper. 
	Then either there exists an uncountable antichain $X \subseteq T$ such that 
	$\meet(X) \subseteq T \setminus K$, or else there is a condition in $\partial \Hcal (K)$ 
	which forces that there is an uncountable antichain $Y \subseteq T$ 
	such that $\meet(Y) \subseteq K$.
\end{lem}

\begin{proof}
	Fix $M \prec H(\chi)$ 
	which is countable such that $T$, $K$, and $\partial \Hcal (K)$ are in $M$, and let 
	$N = M \cap H({2^{\omega_1}}^+)$. 
	Note that $N \in E_1$. 
	Fix any $x \in T_{M \cap \omega_1}$ and define 
	$p = (\{ N \}, \{ x \})$. 
	First, assume that $p$ is not in $\partial \Hcal (K)$. 
	It easily follows that there exists $P \in N \cap E_0$ which rejects $\{ x \}$. 
	Then $P$ rejects the $1$-tuple $(x \restriction (P \cap \omega_1))$, so 
	Lemma \ref{reject iff R} implies that $x \restriction (P \cap \omega_1) \in R$ 
	for some $R \in P \cap \Rscr_1$. 
	It follows that $R$ is uncountable, so by Lemma \ref{R1 uncountable} we are done. 
	Secondly, suppose that $p$ is in $\partial \Hcal (K)$. 
	Let $G$ be a generic filter on $\partial \Hcal (K)$ such that $p \in G$. 
	Since $\partial \Hcal (K)$ is canonically proper, $M[G] \cap V = M$, and 
	therefore $x \notin M[G]$. 
	Since $M[G]$ is an elementary submodel of $H(\chi)^{V[G]}$, 
	by elementarity the set 
	$Y = \bigcup \{ X_q : q \in G \}$ is in $M[G]$, 
	and as $x \in Y \setminus M[G]$, $Y$ is uncountable. 
	By the definition of $\partial \Hcal (K)$, $Y$ is an uncountable antichain and 
	$\meet(Y) \subseteq K$.
\end{proof}

If $\partial \Hcal (K)$ is not canonically proper, then there are two alternatives: 
either (1) $\partial (K)$ is canonically proper 
but $\partial \Hcal (K)$ is not, or (2) both $\partial (K)$ and $\partial \Hcal (K)$ 
are not canonically proper. 
In \cite{linear_basis}, it is shown that both alternatives are false assuming \textsf{PFA}. 
We now connect the core ideas of the proofs of these facts with the idea 
of stability.

Concerning alternative (1), the following lemma was proven as a part of the proof of 
Lemma 5.22 of \cite{linear_basis}.

\begin{lem} \label{partialK not partialHK}
	Suppose that $\partial (K)$ is canonically proper but 
	$\partial \Hcal (K)$ is not. 
	Then there exist $r \in \partial (K)$, $E \in \Ecal$, and 
	$n \in \omega$ positive such that $r$ forces in $\partial (K)$ that there 
	exists a sequence $\Seq{ X_\alpha : \alpha \in \omega_1 }$ of elements 
	of $T^{\otimes n}$ satisfying:
	\begin{enumerate}
		\item \label{none_rejected} for all $P \in E$ and for all $\alpha \in \omega_1$, 
		$P$ does not reject $X_\alpha$;
		\item \label{meet_antichain} for all $\alpha < \beta$ there exists $i < n$ such that 
		$X_\alpha(i) \meet X_\beta(i) \notin K$.
	\end{enumerate} 
\end{lem}

\begin{lem} \label{partialK not stable}
	Suppose that $\partial (K)$ is canonically proper but $\partial \Hcal (K)$ is not. 
	Then $\partial (K)$ is not stable for $T$ and $K$.
\end{lem}

\begin{proof}
	Fix $r$, $E$, and $n$ as in Lemma \ref{partialK not partialHK}, and   
	let $\Seq{ \dot X_\alpha : \alpha \in \omega_1 }$ be a sequence of 
	$\partial (K)$-names for which $r$ 
	forces properties (\ref{none_rejected}) and (\ref{meet_antichain}) of that lemma to hold. 
	Define $S$ as the set of $X \in T^{\otimes n}$ such that for all $P \in E$, 
	$P$ does not reject $X$. 
	Then $r$ forces that for all $\alpha \in \omega_1$, $\dot X_\alpha \in \check S$. 
	In particular, $S$ is uncountable. 
	By Lemma \ref{stationarily many DNR} and (\ref{none_rejected}), $S \in \Rscr_n^{\perp}$, and by 
	Lemma \ref{Rnperp implies two in K} and (\ref{meet_antichain}), $r$ forces that $S \notin \Rscr_n^{\perp}$. 
	So $r$ destabilizes $T$ and $K$.
\end{proof}

Concerning alternative (2), a proof of the following lemma 
appears as a part of the proof of 
Lemma 5.31 of \cite{linear_basis}.

\begin{lem}
	Assume that the Key Lemma holds and 
	$\partial (K)$ is not canonically proper. 
	Then there exist $E \in \Ecal$, $n \in \omega$ positive, 
	and a c.c.c.\ forcing $\q$ such that $\q$ forces that there exists a sequence  
	$\Seq{ X_\alpha : \alpha \in \omega_1 }$ 
	of elements of $T^{\otimes n}$ satisfying:
	\begin{enumerate}
		\item for all $P \in E$ and for all $\alpha \in \omega_1$, 
		$P$ does not reject $X_\alpha$;
		\item for all $\alpha < \beta$ there exists $i < n$ such that 
		$X_\alpha(i) \meet X_\beta(i) \notin K$.
	\end{enumerate} 
\end{lem}

Using Lemma \ref{phiTK implies key lemma}, the next lemma follows by the same 
argument as we gave for Lemma \ref{partialK not stable}.

\begin{prop} \label{ccc not stable}
	Suppose that $\varphi(T,K)$ holds and $\partial (K)$ is not canonically proper. 
	Then there exists a c.c.c.\ forcing which is not stable for $T$ and $K$.
\end{prop}

\section{The Main Results} \label{The Main Results}

We now have all of the tools at our disposal which are 
needed to prove the main theorems of the article.

\begin{thm} \label{main theorem 1}
	Let $T$ be a special Aronszajn tree which is coherent, uniform, and binary. 
	Then for any $K \subseteq T$, 
	there exists a proper forcing which forces that there is 
	an uncountable antichain $A \subseteq T$ such that either $\meet (A) \subseteq K$ 
	or $\meet (A) \subseteq T \setminus K$.
\end{thm}

\begin{proof}
	Applying Theorem \ref{stabilize}, let $\p$ be a proper forcing which forces that every 
	proper forcing is stable for $T$ and $K$. 
	For each positive $n \in \omega$, let $\dot \Fcal_n$ be a $\p$-name for 
	the set $\Rscr_n \cup \Rscr_n^{\perp}$ as computed in $V^{\p}$. 
	By Lemma \ref{forcing multiple psi}, we can fix a $\p$-name $\dot \q$ 
	for a proper forcing which forces $\psi(\dot \Fcal_n)$ for all positive $n \in \omega$. 
	Applying Lemma \ref{forcing phiTK}, let $\dot \Rbb$ be a $\p * \dot \q$-name 
	for a proper forcing which forces $\varphi(T,K)$. 
	Since $\psi$ is upwards absolute, for all positive $n \in \omega$, 
	$\p * \dot \q * \dot \Rbb$ forces $\psi(\Fcal_n)$. 
	By Lemma \ref{preserving stable}, $\p * \dot \q * \dot \Rbb$ forces 
	that every proper forcing is stable for $T$ and $K$. 
	By Propositions \ref{partialK not stable} and \ref{ccc not stable}, 
	the forcing $\p * \dot \q * \dot \Rbb$ forces that $\partial \Hcal (K)$ is 
	canonically proper. 
	Considering the two cases described in Lemma \ref{partialHK done}, 
	let $\dot \Sbb$ be a $\p * \dot \q * \dot \Rbb$-name which is forced 
	to be equal to: (a) in the first case, the trivial forcing, and 
	(b) in the second case, $\partial \Hcal (K)|_{r}$ for some 
	condition $r$ which forces the existence of an uncountable antichain 
	all of whose meets are in $K$. 
	The forcing $\p * \dot \q * \dot \Rbb * \dot \Sbb$ is as required.
\end{proof}

The nontrivial content of the following two lemmas is stated without proof in \cite{AS}; details can be found in 
\cite[Proposition 8.7]{lipschitz}.

\begin{lem} \label{CAT iff CATT}
	Assume $\textsf{MA}_{\omega_1}$ and any two normal Aronszajn trees 
	are club isomorphic. Then the following statements are equivalent:
	\begin{enumerate}
		\item There exists a normal Aronszajn tree $T$ such that $\textsf{CAT}(T)$ holds.
		\item $\textsf{CAT}$.
	\end{enumerate}
\end{lem}

We note that by \cite{AS}, the proper forcing axiom for proper forcings of size at most $\omega_1$ implies 
the assumptions of Lemma \ref{CAT iff CATT}.

By \emph{Baumgartner's axiom} we mean the statement that any two $\omega_1$-dense sets of reals 
are isomorphic (\cite{reals_iso}). 

\begin{lem} \label{equivalences of CAT}
	Assume $\textsf{MA}_{\omega_1}$, 
	any two normal Aronszajn trees are club isomorphic, and Baumgartner's axiom.
	Then the following are equivalent:
	\begin{enumerate}
		\item \textsf{CAT};
		\item Shelah's conjecture;
		\item there exists a five element basis for the uncountable linear orders 
		consisting of $\omega_1$ and its converse $\omega_1^*$, 
		any set of reals of size $\omega_1$, 
		and any Countryman line $C$ and its converse $C^*$.
	\end{enumerate}
\end{lem}

Putting these together, we obtain:

\begin{thm} \label{BPFA implies CAT}
	\textsf{BPFA} implies \textsf{CAT}.
\end{thm}

\begin{proof}
	Let $T$ be a special Aronszajn tree which is coherent, uniform, and binary. 
	For any $K \subseteq T$, the statement that there exists an uncountable antichain 
	whose meets are either all in $K$ or all in $T \setminus K$ is a $\Sigma_1$-statement 
	with parameters in $H(\omega_2)$. 
	By Theorem \ref{main theorem 1}, there exists a proper forcing which forces 
	this statement, so by \textsf{BPFA}, this statement is true. 
	As \textsf{BPFA} implies the assumptions of Lemma \ref{CAT iff CATT}, we are done.
\end{proof}

Since the assumptions of Lemma \ref{equivalences of CAT} follow from \textsf{BPFA}, 
we immediately obtain the following corollary.

\begin{cor}
	\textsf{BPFA} implies Shelah's conjecture and the existence of a five element 
	basis for the uncountable linear orders.
\end{cor}

\begin{thm} \label{inaccessible to five}
	If $\kappa$ is inaccessible, then in $L$ there exists 
	a proper forcing which forces the proper forcing axiom for proper 
	forcings of size at most $\omega_1$, 
	Baumgartner's axiom, Shelah's conjecture, and $\varphi$.
\end{thm}

\begin{proof}
	Let $\kappa$ be an inaccessible cardinal. 
	Then $\kappa$ is inaccessible in $L$ and $\Diamond_\kappa$ holds. 
	We define a countable support forcing iteration $\p$ of length $\kappa$ of 
    forcings of size less than $\kappa$. 
	In order to arrange that $\p$ forces $\varphi$, we repeat the argument of 
	Theorem \ref{inaccessible equiconsistent to phi}. 
	Since $\kappa$ remains inaccessible after forcing with any initial segment of $\p$, 
    and hence $V_\kappa$ is a model of \textsf{ZFC} in any generic extension 
    by such an initial segment, 
	we can apply Theorem \ref{main theorem 1} to force instances of $\textsf{CAT}$ 
    with forcings of size less than $\kappa$, 
	bookkeeping to make sure all $T$ and $K$ which appear in $V^\p$ are handled. 
	Similarly, any required instance of Baumgartner's axiom or the proper 
	forcing axiom for proper 
	forcings of size at most $\omega_1$ can be handled by standard bookkeeping. 
	Now we are done by Lemma \ref{equivalences of CAT}.
\end{proof}

\begin{cor} \label{equiconsistency}
	The following are equiconsistent:
	\begin{enumerate}
		\item There exists an inaccessible cardinal.
		\item There exists a five element basis for the uncountable linear orders and $\varphi$ holds.
	\end{enumerate}
\end{cor}

\section{Open Problems}

We close the article by stating some questions related to Shelah's conjecture which remain open.

\begin{question}[\cite{linear_basis}] \label{q1}
	Does Shelah's conjecture have any large cardinal consistency strength?
\end{question}

It is not known, for example, if the existence of a Kurepa tree implies the 
failure of Shelah's conjecture. 
If it does, then Corollary \ref{equiconsistency} provides the optimal consistency 
strength result.

\begin{question} \label{q2}
	Does Shelah's conjecture imply the non-existence of a Kurepa tree, 
	or Aronszajn tree saturation?
\end{question}

Similar questions can be asked by replacing Shelah's conjecture with 
any of the statements listed in the conclusion of Lemma \ref{equivalences of CAT}.

\textsf{BPFA} implies the non-existence of a Kurepa tree, 
but we do not know if it implies the related statement of Aronszajn tree saturation.
Under \textsf{BPFA}, $\varphi$ is equivalent to Aronszajn tree saturation.

\begin{question}[\cite{con_linear_basis}] \label{q3}
	Does \textsf{BPFA} imply Aronszajn tree saturation?
\end{question}

\section{Acknowledgements}

The second author acknowledges support from the National Science Foundation under 
Grants DMS-1854367 and DMS-2153975. 
Part of the work for this article was conducted while the first author was visiting 
the second author at Cornell University in December 2024. 
The first author thanks the second author for funding his visit and for his 
hospitality. 
The first author also thanks Šárka Stejskalová for discussions about the idea of a subtree base 
which motivated some of the methods developed in this article.

\end{document}